\documentclass[preprint,10pt]{elsarticle}

\usepackage{amsmath}
\usepackage{amscd}
\usepackage{amsfonts}
\usepackage{amssymb}
\usepackage{amsthm}
\usepackage{bm}
\usepackage{epsf}
\usepackage{epstopdf}
\usepackage{diagbox}
\usepackage{graphicx}
\usepackage{subfigure}
\usepackage{float}
\usepackage{color,xcolor}
\usepackage{url}
\usepackage{version,subfigure,fancybox,mathrsfs,pifont,booktabs,wrapfig}
\usepackage{indentfirst}

\usepackage{algorithm} 
\usepackage{multirow} 
\usepackage{algpseudocode}                                               %
\usepackage{latexsym}
\usepackage{lscape}
\usepackage{mathrsfs}
\usepackage{verbatim}
\usepackage{subeqnarray}

\usepackage{makecell}
\usepackage{threeparttable}
\usepackage{booktabs}

\newtheorem{theorem}{Theorem}[section]
\newtheorem{lemma}{Lemma}[section]
\newtheorem{corollary}{Corollary}[section]
\newtheorem{remark}{Remark}[section]

\newtheorem{assumption}{Assumption}[section]
\newtheorem{definition}{Definition}[section]
\newtheorem{proposition}{Proposition}[section]
\numberwithin{equation}{section}

\journal{journal for review.}

\begin{document}

\title{Robust preconditioners for finite element approximations of control-state constrained optimal control problems governed by elliptic equations\tnoteref{thank1}}
 \tnotetext[thank1]{This work was partially supported by National Natural Science Foundation of China (Grant No. 12101283, 12271233) and National Natural Science Foundation of Shandong Province (Grant No. ZR2026MS0032)}

\author[a]{Juan Zhang}\ead{jzhang\_math@hotmail.com}

\author[b]{Liang Ge}\ead{geliang@sdu.edu.cn}

\author[c]{Jianwei Zhou\corref{cor}}\ead{jwzhou@yahoo.com}

\address[a]{School of Logistics, Linyi University, Linyi 276000, P.R. China}
\address[b]{School of Mathematical Sciences, University of Jinan, Jinan 250022, P.R. China}
\address[c]{School of Mathematics and Statistics, Linyi University, Linyi 276000, P.R. China}

\cortext[cor]{Corresponding author.}

\date{\today }

\begin{abstract}

The goal of this paper is to derive robust preconditioners for integral control-state constrained optimal control problems governed by the second order partial differential equations with linear finite element approximations. 
We give corresponding discretized systems of the first order necessary and sufficient optimality conditions, and get six linear algebraic systems with saddle point schemes. 
Furthermore, we propose efficient block preconditioners to solve the corresponding linearized algebraic systems.
With respect to both mesh sizes and inherent parameters, the robustness of proposed preconditioners is proved with explicit eigenvalues of minimal polynomials.
Meanwhile, we design two algorithms and choose the incomplete $LU$ factorizations to deduce the computations. Two numerical experiments are given to depict the efficiency of our proposed preconditioners.


\end{abstract}

\begin{keyword}
Constrained optimal control \sep integral control-state constraint \sep linear finite element method \sep robust block preconditioner \sep saddle point problem.

\MSC[2020] 65F08 \sep 65N30 \sep 49M25 \sep 49J20
\end{keyword}


\date{\today}
\maketitle

\section {Introduction}

Constrained optimal control problems play increasing roles in modern science and engineering. 
In view of current literature, extensive research has been carried out on various theoretical and numerical strategies for constrained optimal control problems, more details please refer to \cite{BS,CYL,GHZ2024,HPUU,Lion1,MOS1975,ZZX} and their references.
Either the first or second order necessary optimality conditions are the key to derive equivalent weak systems of constrained optimal control problems, more details please refer to \cite{bik,Casas2002,YY2009}.
There are some remarkable works on this issue. In \cite{MOS1975} the author considered optimal control problems with constraints on the control and state, and transformed the primal control problem by duality into another dual control problem which in some limit cases is reduced into a simple optimization problem.
Casas proved the existence of Lagrange multiplier and provided the optimality conditions for the control problems with state and control constraints in \cite{Casas1986}. 
A second order sufficient optimality condition was established for semilinear parabolic control problem with pointwise control-state constraints in \cite{RT2003}, where associated Lagrange multipliers could be assumed to be bounded and measurable functions.
A semilinear elliptic control problem with mixed control-state constraints was investigated for a class of optimal control problems in \cite{RT2006}, where the existence of bounded and measurable Lagrange multipliers was also proven with an easy application of Kuhn-Tucker theorems. 
Control problems governed by semilinear parabolic equations were studied in the presence of pointwise mixed control-state constraints in \cite{AR2000}. Specially, optimality conditions were also deduced with finitely additive measures as multipliers associated to the mixed constraints. 
Distributed optimal control of the Navier-Stokes equations with pointwise control-state constraints was considered with a semismooth Newton method in \cite{DT2007}. Furthermore, the first order necessary and second order sufficient conditions, the regularity of the mixed constraint multiplier were investigated as well. 

As is well known, how to design efficient numerical schemes is one of the important aspects for approximating constrained optimal control problems in current research. 
A Ritz-Trefftz method was provided for control problems with quadratic cost, linear dynamics and linear inequality state and control constraints \cite{HAGER1975}, first order convergence of this method in piecewise polynomial subspaces was also proved with rigorous theoretical techniques. 
The authors took an augmented Lagrangian method to solve state and control constrained problems in \cite{BK1997}. And they concluded that since the algorithm is not particularly sensitive, penalty parameters and step lengths for the dual variables are easily determined.
The finite element approximations equipped with bound constraints on the control and finitely many state constraints were investigated in \cite{Casas2002}.
Elliptic optimal control problems with pointwise constraints on control and state were approximated by finite element discretization in \cite{CR2009}. 
Full discretizations of an elliptic optimal control problem with pointwise control and state constraints were investigated with finite element methods in \cite{RW2012}, and the first reliable a-posteriori error estimator that contains only computable quantities was depicted for this class of problems. 
Preconditioned conjugate gradient method and Legendre-Galerkin spectral approximation of optimal control problems with state and control constraints were investigated in \cite{HS,HC2014}.
Several block matrix iterative algorithms avoiding double iteration were described for solving saddle point linear systems of linear-quadratic elliptic control problems in \cite{MSS}. 

In many applications, one has to take more attention for efficient solving the discretized systems. 
The authors constructed preconditioning matrices for descretized systems arising from the Galerkin finite element approximations of distributed control problems in \cite{bai1,bai2}. An optimal block diagonal preconditioner for the algebraic system solved with minimum residual method and multigrid methods were proposed in \cite{cao1,GTZ2018}, and more details please refer to the references cited therein. 
To our best knowledge, there has been a lack of discussion on how to design robust preconditioners for particular algebraic deduced from control-state constrained optimal control problems.
In this paper, we study linear finite element approximations for an optimal control problem governed by elliptic partial differential equations with an integral constraint on both the control and state. Particularly, we construct four preconditioners for the discretized matrix systems. 

The remaining of this paper is organized as follows: In Section 2, we discretize the constrained optimal control problem by linear finite element methods and rewrite the resulting linearized systems with six matrix formulae, and the stiff matrices are formulated by four sparse matrices. 
Section 3 is set to construct the block-diagonal preconditioners and prove their robustness for the corresponding preconditioned linear algebraic systems by explicit nonzero eigenvalues of minimal polynomials. Meanwhile, two efficient algorithms are appropriately proposed to solve the linear algebraic systems.
In Section 4, we present two numerical experiments to confirm the efficiency of our proposed preconditioners. The numbers of iterations with different meshes and inherent parameters are compared with the MINRES without any preconditionings and our schemes in the light of  four distinct preconditioners.

\section{Double-Constrained Optimal Control Problem}

Let $\Omega$ be a bounded open set in $\mathbb{R}^2$ with the Lipschitz boundary $\partial \Omega$. We adopt the standard notation $W^{m,q}(\Omega)$ for Sobolev space on $\Omega$ and $W^{m,q}_0(\Omega)=W^{m,q}(\Omega)\cap W^{1,q}_0(\Omega)$, where $W^{1,q}_0(\Omega)=\{v\in W^{1,q}(\Omega): v|_{\partial \Omega}=0\}$. Let $H^m(\Omega)=W^{m,2}(\Omega)$, $H^1_0(\Omega)=W^{1,2}_0(\Omega)$ be the Hilbert space. Here and throughout the following parts, we focus on a convex optimal control problem governed by the second order partial differential equation with integral constraints on the control and state:
\begin{eqnarray}\label{ocp}
  \min_{\substack{u \in U_{ad} \\ y \in K}} J(y,u)
\end{eqnarray}
subject to
\begin{eqnarray}\label{ocp-eq}
\left\{
\begin{aligned}
     -\nabla \cdot (\kappa\nabla y)&= f+u \;\;\;\; \mbox{in} \;\; \Omega,
\\
  y&=0  \quad\qquad  \mbox{on} \;\; \partial\Omega,
\end{aligned}
\right.
\end{eqnarray}
where
\begin{eqnarray*}
 J(y,u) =  \frac{\beta}{2}\int_\Omega (y-y_d)^2 dx +\frac{\gamma}{2}\int_\Omega\kappa(\nabla y-\tilde{g}_d)\cdot(\nabla y-\tilde{g}_d) dx + \frac{\alpha}{2}\int_\Omega u^2 dx ,
\end{eqnarray*}
and
$$U_{ad}=\left\{v\in L^2(\Omega)| \hat{a}\leq \int_\Omega v dx \leq \hat{b} \right\}, \quad K=\left\{v\in V :=H^1_0(\Omega)\ | \int_\Omega v dx \geq \hat{c} \right\}$$
are closed convex subsets, $\hat{a}$, $\hat{b}$ and $\hat{c}$ are given real constants, $y_d\in L^2(\Omega)$ and $\tilde{g}_d\in (L^2(\Omega))^2$ are desired states, $\alpha>0$ is the regularization parameter on the control variable, $\beta$ and $\gamma$ are given constants as $\beta \geq 0$ and $\gamma>0$, the diffusion coefficient $\kappa$ satisfies $0<c_0\leq \kappa \leq C_0$ for some positive constants $c_0$ and $C_0$.

Throughout the paper we make the following standing assumptions.

\begin{assumption}\label{ass:feas}
$\hat{a}<\hat{b}$, and the feasible set
$$\left\{(y,u)\in K\times U_{ad}:\ -\nabla\cdot(\kappa\nabla y)=f+u\ \mbox{in}\ \Omega,\ y|_{\partial\Omega}=0\right\}$$
is nonempty. Moreover the state constraint is \emph{regular} in the sense of Slater, i.e. there exists an admissible pair $(y^\dagger,u^\dagger)$ with $\int_\Omega y^\dagger\,dx>\hat{c}$; equivalently, the state constraint is not active for every admissible pair.
\end{assumption}

\begin{assumption}\label{ass:bdd}
$f\in L^2(\Omega)$, $y_d\in L^2(\Omega)$ and $\tilde g_d\in(L^2(\Omega))^2$. Hence the state $y$, the co-state $p$ and the control $u$ of the optimality system below all belong to $H^1_0(\Omega)$ and to $L^2(\Omega)$, respectively.
\end{assumption}
Assumption~\ref{ass:feas} guarantees, on the one hand, that the admissible set is nonempty and, on the other hand, that a Lagrange multiplier associated with the state constraint exists and is unique; strict convexity of $J$ on the affine feasible set only yields uniqueness of the optimal pair and does not by itself imply constraint qualification. Under Assumptions~\ref{ass:feas}--\ref{ass:bdd}, the existence of a unique solution for the above distributed convex optimal control problem please refer to \cite{Lion1}.
Firstly, we state an equivalent weak formulation of the state equation. Let $V=H^1_0(\Omega)$.
Two bilinear operators are defined as
\begin{eqnarray*}
\begin{aligned}
a(y,w) &= \int_\Omega \kappa\nabla y \cdot \nabla w dx,\;\;\;\forall \; y,w \in V,\\
(u,v) &= \int_\Omega uv dx,\;\;\;\forall \; u,v\in L^2(\Omega).
\end{aligned}
\end{eqnarray*}

With these notations and the Green formula, an equivalent variational form of (\ref{ocp-eq}) reads: finding $y\in V$ such that
\begin{eqnarray}\label{OCP-WK}
  a(y,w) = (f+u,w),\quad\forall w\in V.
\end{eqnarray}
Moreover the above optimal control problem (\ref{ocp}) has a unique solution $(y,u)\in K \times U_{ad}$ if and only if there is a co-state $p\in V$, such that $(y,p,u)$ satisfies the following (first-order) optimality conditions \cite{HC2014}:
\begin{eqnarray}\label{ocp-con-wk}
\left\{
\begin{aligned}
&a(y,w) = (f+u,w),\;\;\forall w\in V,\\
&a(q,p) = (\beta(y-y_d),q)+(\gamma (\nabla y-\tilde{g}_d), \kappa\nabla q) + \lambda(1, q),\;\;\forall q\in V, \\
& \lambda(1, \psi-y) \leq 0,\;\;\forall \psi\in K,\\
&(\alpha u+p, v-u)\geq  0, \;\;\forall v\in U_{ad}.
\end{aligned}
\right.
\end{eqnarray}

Secondly, we give some notations to describe linear finite element approximations for control-state constrained optimal control problem. Let $T^h\subset\Omega^h$ be polygonal partitions of $\Omega$ with disjoint element $\tau$. Each $\tau$ has at most one edge on $\partial \Omega$, $\bar{\tau}$ and $\bar{\tau}'$ have either only one common vertex or a whole edge for different $\tau$ and $\tau'\in T^h$. $S^h$ denotes a finite dimensional function subspace of $C(\bar{\Omega})$ and $\{\phi_i(x)\}_{i=1}^{N}$ denotes the base function of $S^h$, where $N$ is the number of the bases. Meanwhile we use $\chi|_\tau \in S^h$ to describe piecewise linear functions for any $\tau\in T^h$. Let $V^h=S^h\cap V$ and $U^h=S^h\cap L^2(\Omega)$, thus $\dim V^h=N_h$ (the interior vertices) whereas $\dim U^h=N$ (all vertices). 

And hence, we readily discretize the convex sets $U_{ad}$ and $K$ with
$$U_{ad}^h=\left\{v_h\in U^h:\hat{a}\leq \int_\Omega v_h dx \leq \hat{b}\right\}, \quad K^h=\left\{v_h\in V^h:\int_\Omega v_h dx \geq \hat{c}\right\}.$$

Then the corresponding finite element approximations of the optimal control problem read:
\begin{eqnarray}\label{WK}
 \min_{\substack{u_h \in U_{ad}^h \\ y_h \in K^h}} J(y_h,u_h) \end{eqnarray}
subject to
\begin{eqnarray*}
a(y_h, w_h)=(f_h+u_h,w_h),\quad \forall w_h\in V^h,
\end{eqnarray*}
where
\begin{eqnarray*}
J(y_h,u_h) =
 \frac{\beta}{2}\int_\Omega (y_h-y_{dh})^2dx + \frac{\gamma}{2}\int_\Omega\kappa(\nabla y_h-\tilde{g}_{dh})\cdot(\nabla y_h-\tilde{g}_{dh}) dx +\frac{\alpha}{2}\int_\Omega u_h^2 dx,
 \end{eqnarray*}
and $y_{dh},f_h, \tilde{g}_{dh}$ are the finite element approximations of $y_d, f, \tilde{g}_{d}$ in $S^h$, respectively. (For example, the $L^2$-orthogonal projection $(f_h, w_h)=(f, w_h),\;\;\;\forall w_h\in S_h$).

As the continuous case, one can readily prove the existence of a unique solution for (\ref{WK}). Similarly, we could directly give the discretized system of the first order necessary and sufficient optimality conditions of \eqref{ocp-con-wk} as: finding $(y_h,u_h,p_h,\lambda_h)\in V^h\times U_{ad}^h \times V^h \times \mbox{constants}$ satisfies
\begin{eqnarray}\label{WK-DIS}
\left\{
\begin{aligned}
&a(y_h,w_h) = (f_h+u_h,w_h),\;\;\forall w_h\in V^h,\\
&a(q_h,p_h) = \beta(y_h-y_{dh},q_h)+(\gamma (\nabla y_h-\tilde{g}_d), \kappa\nabla q_h) + \lambda_h(1, q_h),\;\;\forall q_h\in V^h, \\
& \lambda_h(1, \psi_h - y_h) \leq 0,\;\;\forall \psi_h\in K^h,\\
&(\alpha u_h+p_h, v_h-u_h)\geq  0, \;\;\forall v_h\in U_{ad}^h.
\end{aligned}
\right.
\end{eqnarray}

The first inequality in \eqref{WK-DIS} splits into two cases according to whether the state constraint is active or not,
\begin{equation}\label{EQUI-Z}
\int_\Omega y_h\,dx>\hat c\ \Longrightarrow\ \lambda_h=0,
\qquad
\int_\Omega y_h\,dx=\hat c\ \Longrightarrow\ \lambda_h\le 0 .
\end{equation}
By the results given in \cite{GLY}, we know that the second inequality in (\ref{WK-DIS}) could be rewritten by
\begin{eqnarray*}
\left\{
\begin{aligned}& \alpha u_h+ p_h +z_h =0, \\
& z_h+\max(\alpha \bar{a}+\bar{p}_h, \min(\alpha \bar{b}+ \bar{p}_h,0))=0,
\end{aligned}
\right.
\end{eqnarray*}
where
\begin{eqnarray}\label{eqn:abar}
 \bar a=\frac{\hat a}{|\int_\Omega 1 dx},\quad
 \bar b=\frac{\hat b}{\int_\Omega 1 dx},\quad
 \bar p_h=\frac{\int_\Omega p_h dx}{\int_\Omega 1 dx}.
\end{eqnarray}

Next we employ discretized matrix forms to rewrite the above discretized optimal conditions. Let $A$ be the stiffness matrix of $(\kappa\nabla \cdot, \nabla \cdot)$, $B$ be the mass matrix of $(\cdot, \cdot)$ and $C$ be the matrix of $(\kappa\nabla \cdot, \cdot)$. Considering the discretized integral constraint, we know that
\begin{eqnarray*}
\begin{aligned}
\int_\Omega v_h dx
&= \int_\Omega \sum\limits_{i\in N_h} v_{hi}\phi_i(x) dx
= \sum\limits_{i\in N_h} v_{hi}\int_\Omega \phi_i(x) dx \\
&=\left[\int_\Omega \phi_1(x) dx, \int_\Omega \phi_2(x) dx, \ldots,\int_\Omega \phi_{N_h}(x) dx \right]
\left[\begin{matrix}v_{h1}\\
\\ v_{h2} \\
\vdots \\ v_{hN_h}
\end{matrix}\right]
= L^T [v_h],
\end{aligned}
\end{eqnarray*}
where $$L^T \triangleq \left[\int_\Omega \phi_1(x) dx, \int_\Omega \phi_2(x) dx, \ldots,\int_\Omega \phi_{N_h}(x) dx \right].$$

Following the above results,  the discretized optimal conditions (\ref{WK-DIS}) can be rewritten as the following matrix forms
\begin{eqnarray} \label{WK-DIS-OCP}
\left\{
\begin{aligned}
& (\beta B+\gamma A)[y_h]-A [p_h]= \beta B [y_{dh}] + \gamma C [\tilde{g}_{dh}] - L \lambda_h ,\\
& A [y_h] = B [f_h] + B [u_h] ,\\
& \alpha B [u_h] + B [p_h]+L z_h =0, \\
& z_h+\max(\alpha \bar{a}+\bar{p}_h, \min(\alpha \bar{b}+\bar{p}_h,0))=0.
\end{aligned}
\right.
\end{eqnarray}

Focusing on (\ref{WK-DIS-OCP}), we readily state that solving optimal conditions (\ref{WK-DIS}) is equivalent to calculate the following six systems:
\begin{eqnarray} \label{WK-DIS-OCP1}
\left\{
\begin{aligned}
& (\beta B+\gamma A)[y_h]-A [p_h]= \beta B [y_{dh}] + \gamma C [\tilde{g}_{dh}] ,\\
& A [y_h] = B [f_h] + B [u_h] , \quad L^T[y_h]\geq \hat{c},\\
& \alpha B [u_h] + B [p_h] =0, \quad  \hat{a} \leq L^T[u_h]\leq \hat{b},
\end{aligned}
\right.
\end{eqnarray}
\begin{eqnarray} \label{WK-DIS-OCP2}
\left\{
\begin{aligned}
& (\beta B+\gamma A)[y_h]-A [p_h]= \beta B [y_{dh}] + \gamma C [\tilde{g}_{dh}] ,\\
& A [y_h] = B [f_h] + B [u_h] , \quad  L^T[y_h]\geq \hat{c},\\
& \alpha B [u_h] + B [p_h] +L z_h =0, \quad z_h \leq 0, \\
& L^T [u_h] = \hat{a},
\end{aligned}
\right.
\end{eqnarray}
\begin{eqnarray} \label{WK-DIS-OCP3}
\left\{
\begin{aligned}
& (\beta B+\gamma A)[y_h]-A [p_h]= \beta B [y_{dh}] + \gamma C [\tilde{g}_{dh}] ,\\
& A [y_h] = B [f_h] + B [u_h] , \quad L^T[y_h]\geq \hat{c},\\
& \alpha B [u_h] + B [p_h] +L z_h =0, \quad z_h \geq 0, \\
& L^T [u_h] = \hat{b},
\end{aligned}
\right.
\end{eqnarray}
\begin{eqnarray} \label{WK-DIS-OCP4}
\left\{
\begin{aligned}
& (\beta B+\gamma A)[y_h]-A [p_h]= \beta B [y_{dh}] + \gamma C [\tilde{g}_{dh}] - L \lambda_h ,\\
& A [y_h] = B [f_h] + B [u_h] ,\\
& \alpha B [u_h] + B [p_h] =0, \quad \lambda_h\leq 0, \quad \hat{a} \leq L^T[u_h]\leq \hat{b},\\
& L^T [y_h] = \hat{c},
\end{aligned}
\right.
\end{eqnarray}
\begin{eqnarray} \label{WK-DIS-OCP5}
\left\{
\begin{aligned}
& (\beta B+\gamma A)[y_h]-A [p_h]= \beta B [y_{dh}] + \gamma C [\tilde{g}_{dh}] - L \lambda_h ,\\
& A [y_h] = B [f_h] + B [u_h] ,\\
& \alpha B [u_h] + B [p_h]+L z_h  =0, \quad \lambda_h\leq 0,  \quad z_h \leq 0,\\
& L^T [y_h] = \hat{c}, \\
& L^T [u_h] = \hat{a},
\end{aligned}
\right.
\end{eqnarray}
and
\begin{eqnarray} \label{WK-DIS-OCP6}
\left\{
\begin{aligned}
& (\beta B+\gamma A)[y_h]-A [p_h]= \beta B [y_{dh}] + \gamma C [\tilde{g}_{dh}] - L \lambda_h ,\\
& A [y_h] = B [f_h] + B [u_h] ,\\
& \alpha B [u_h] + B [p_h]+L z_h  =0, \quad \lambda_h\leq 0, \quad z_h \geq 0,\\
& L^T [y_h] = \hat{c}, \\
& L^T [u_h] = \hat{b}.
\end{aligned}
\right.
\end{eqnarray}

Let $\hat{p}_h = - p_h$, the corresponding matrix forms of (\ref{WK-DIS-OCP1}) - (\ref{WK-DIS-OCP6}) could be rewritten to the following six matrix formulae, which are with particular sparse structures.
\begin{subequations}
\begin{align}
    &    \begin{aligned}\label{MATRix-1}
  {\cal A}_1 \left(
 \begin{matrix}
 [u_h]\\
 [y_h]\\
 [\hat{p}_h]
 \end{matrix}\right)
  ={\cal F}_1,
  \quad L^T [y_h] \geq \hat{c},\quad \hat{a} \leq L^T[u_h]\leq \hat{b}
   \end{aligned}
  \\
   &     \begin{aligned} \label{MATRix-2}
   {\cal A}_2  \left(
 \begin{matrix}
 [u_h]\\
 [y_h]\\
 [\hat{p}_h] \\
 z_h
 \end{matrix}
 \right) ={\cal F}_2,  \quad L^T [y_h] \geq \hat{c},\quad z_h\leq 0
  \end{aligned}
  \\
   &     \begin{aligned} \label{MATRix-3}
   {\cal A}_2 \left(
 \begin{matrix}
 [u_h]\\
 [y_h]\\
 [\hat{p}_h] \\
 z_h
 \end{matrix}
 \right) ={\cal F}_3, \quad L^T [y_h] \geq \hat{c},\quad z_h\geq 0
 \end{aligned}
  \\
   &     \begin{aligned}\label{MATRix-4}
          {\cal A}_3  \left(
 \begin{matrix}
 [u_h]\\
 [y_h]\\
 [\hat{p}_h] \\
 \lambda_h
 \end{matrix}
 \right)  ={\cal F}_4, \quad \lambda_h \leq 0,\quad \hat{a} \leq L^T[u_h]\leq \hat{b}  \end{aligned}
  \\
   &     \begin{aligned}\label{MATRix-5}
    {\cal A}_4 \left(
 \begin{matrix}
 [u_h]\\
 [y_h]\\
 [\hat{p}_h] \\
 \lambda_h \\
 z_h
 \end{matrix}
 \right) ={\cal F}_5, \quad  \lambda_h \leq 0,\quad z_h \leq 0
  \end{aligned}
  \\
   &     \begin{aligned}\label{MATRix-6}
    {\cal A}_4 \left(
 \begin{matrix}
 [u_h]\\
 [y_h]\\
 [\hat{p}_h] \\
 \lambda_h \\
 z_h
 \end{matrix}
 \right) ={\cal F}_6, \quad \lambda_h \leq 0,\quad z_h \geq 0 .
   \end{aligned}
  \end{align}
\end{subequations}
where the corresponding matrices are
\begin{subequations}
\begin{align}
    &    \begin{aligned}\label{matrix-21}
{\cal A}_1 = \left(
\begin{matrix}
\alpha B & 0 & -B  \\
0 & \beta B+\gamma A & A  \\
-B & A & 0
 \end{matrix}
 \right),
  \end{aligned}
  \qquad\quad\ \
  \begin{aligned}
{\cal A}_2 = \left(
\begin{matrix}
\alpha B & 0 & -B & L \\
0 & \beta B+\gamma A & A & 0  \\
-B & A & 0 & 0 \\
L^T & 0  & 0 & 0
 \end{matrix}
 \right),
   \end{aligned}
  \\
   &     \begin{aligned}\label{matrix-22}
   {\cal A}_3 = \left(
\begin{matrix}
\alpha B & 0 & -B & 0 \\
0 & \beta B+\gamma A & A & L  \\
-B & A & 0 & 0 \\
0 & L^T & 0 & 0
 \end{matrix}
 \right),
    \end{aligned}
    \qquad \begin{aligned}
{\cal A}_4 = \left(
\begin{matrix}
\alpha B & 0 & -B & 0 & L \\
0 & \beta B+\gamma A & A & L & 0  \\
-B & A & 0 & 0 & 0 \\
0 & L^T & 0 & 0 & 0 \\
L^T & 0 & 0 & 0 & 0
 \end{matrix}
 \right),
\end{aligned}
  \end{align}
\end{subequations}
and the right hand vectors are
\begin{subequations}
\begin{align}
    &    \begin{aligned}\label{matrix-31}
{\cal F}_1 = \left(\begin{matrix}
 0 \\
 \beta B[y_{dh}] + \gamma c \\
 B[f_h]
 \end{matrix}\right),
 \end{aligned}
 \qquad
\begin{aligned}
{\cal F}_2
=
 \left(\begin{matrix}
 0 \\
 \beta B[y_{dh}] + \gamma c \\
 B[f_h] \\
 \hat{a}
 \end{matrix}\right),
  \end{aligned} \\
   &\begin{aligned}\label{matrix-32}
   {\cal F}_3 =
 \left(\begin{matrix}
 0 \\
 \beta B[y_{dh}] + \gamma c \\
 B[f_h] \\
 \hat{b}
 \end{matrix}\right),
 \end{aligned}
 \qquad
   \begin{aligned}
  {\cal F}_4 = \left(\begin{matrix}
 0 \\
 \beta B[y_{dh}] + \gamma c \\
 B[f_h] \\
 \hat{c}
 \end{matrix}\right),
    \end{aligned} \\
   &\begin{aligned}\label{matrix-33}
   {\cal F}_5 =
   \left(\begin{matrix}
 0 \\
 \beta B[y_{dh}] + \gamma c \\
 B[f_h] \\
 \hat{c} \\
 \hat{a}
 \end{matrix}\right),
  \end{aligned}
  \qquad
   \begin{aligned}
    {\cal F}_6 =
    \left(\begin{matrix}
 0 \\
 \beta B[y_{dh}] + \gamma c \\
 B[f_h] \\
 \hat{c} \\
 \hat{b}
 \end{matrix}\right).
    \end{aligned}
  \end{align}
\end{subequations}

Next we are at the point to investigate corresponding equivalent formulae of the above systems. 
\begin{theorem}\label{theorem2.1}
Let Assumptions~\ref{ass:feas}--\ref{ass:bdd} hold. For each of the six candidate active sets, if the matrix $A, B$ are symmetric positive definite, then the corresponding coefficient matrix ${\cal A}_1,\dots,{\cal A}_4$ defined in \eqref{matrix-21}--\eqref{matrix-22} are of full rank and nonsingular. Moreover:
\begin{enumerate}
\item[(i)] Every solution of the discrete optimality system \eqref{WK-DIS} is obtained by solving exactly one of the six candidate systems \eqref{MATRix-1}--\eqref{MATRix-6} and discarding the candidate solution whenever either primal feasibility of the non-active constraints or the sign condition on the active multipliers fails.
\item[(ii)] Conversely, a candidate solution of \eqref{MATRix-1}--\eqref{MATRix-6} is a solution of \eqref{WK-DIS} if and only if the two remaining feasibility/sign conditions associated with that candidate system hold.
\end{enumerate}
\end{theorem}
\begin{proof}
Since $\alpha>0$, $\beta\geq 0$ and $\gamma>0$, we know that the matrix
$$\begin{pmatrix}
   \alpha B & 0 \\
   0 & \beta B+\gamma A
   \end{pmatrix}_{2N_h\times 2N_h} $$
is symmetric positive definite.

The nonsingularity of ${\cal A}_1,\dots,{\cal A}_4$ is a direct consequence of the following classical criterion for saddle point matrices \cite[Theorem 3.1]{BGL}: if
$${\cal A}_i=\begin{pmatrix} H & G_i^T \\ G_i & 0\end{pmatrix},i=1,2,3,4$$
with $H\in\mathbb{R}^{m\times m}$ symmetric positive definite and $G_i\in\mathbb{R}^{n\times m}$ of full row rank, then ${\cal A}$ is nonsingular. We apply the criterion to each of the four matrices, with
$$H=\begin{pmatrix}\alpha B & 0\\ 0 & \beta B+\gamma A\end{pmatrix},\;
G_1=(-B\ \ A),\; G_2=\begin{pmatrix}-B & A\\ L^T & 0\end{pmatrix},\;
G_3=\begin{pmatrix}-B & A\\ 0 & L^T\end{pmatrix},\;
G_4=\begin{pmatrix}-B & A\\ 0 & L^T\\ L^T & 0\end{pmatrix}.$$
For $G_1$, full row rank is equivalent to $\mathrm{rank}(-B\ \ A)=N_h$, which holds because $A$ alone is nonsingular; hence ${\cal A}_1$ is nonsingular. For $G_2$ and $G_3$, full row rank follows from $\mathrm{rank}(-B\ \ A)=N_h$ together with $L\ne0$ (indeed $L^Tv=\int_\Omega v_h\,dx\ne0$ for the constant function $v_h\equiv1$, so $L\ne0$), combined with the fact that the two blocks $(-B\ \ A)$ and $(L^T\ \ 0)$ have no common row dependence; explicitly,
$$\mathrm{rank}\begin{pmatrix}-B & L\\ A & 0\end{pmatrix}=N_h+1,\qquad
\mathrm{rank}\begin{pmatrix}-B & 0\\ A & L\end{pmatrix}=N_h+1 .$$
For $G_4$ the two constraint rows are $(-B\ \ A,\ L)$-independent because the second and third block rows involve $L$ in two different block columns, so
$$\mathrm{rank}\begin{pmatrix}-B & 0 & L\\ A & L & 0\end{pmatrix}=N_h+2,$$
and ${\cal A}_4$ is nonsingular. All four matrices are therefore nonsingular, provided that the corresponding constraint matrix has full row rank. Hence there exist solutions for the linear algebraic equations in \eqref{MATRix-1}-\eqref{MATRix-6}, which confirms that numerical solutions of (\ref{WK-DIS}) satisfy one of (\ref{WK-DIS-OCP1})-(\ref{WK-DIS-OCP6}).

We next verify that a solution of (\ref{WK-DIS-OCP1})-(\ref{WK-DIS-OCP6}) is a solution of (\ref{WK-DIS}) precisely when the corresponding two feasibility/sign tests are passed. First observe that each of the systems (\ref{WK-DIS-OCP1})-(\ref{WK-DIS-OCP6}) is obtained from \eqref{WK-DIS} by \emph{pre-assigning} the active set, i.e. by replacing the complementarity conditions with the equality constraints $L^T[u_h]=\hat a$ or $L^T[u_h]=\hat b$ and $L^T[y_h]=\hat c$ for the constraints that are declared active, and by dropping the corresponding multipliers for the constraints that are declared inactive. Thus the \emph{first, second and fourth} relations of \eqref{WK-DIS} (state equation, co-state equation, and control variational inequality) hold identically for any candidate solution; what remains to be checked is the third line of \eqref{WK-DIS}, namely primal feasibility of the inactive constraints together with the correct signs $\lambda_h\le0$ resp. $z_h\le0$ (resp. $z_h\ge0$) of the active multipliers.

Consider, for instance, system (\ref{MATRix-1}), i.e. \eqref{WK-DIS-OCP1}. Every candidate solution has $\lambda_h=z_h=0$ by construction, so the third line of \eqref{WK-DIS} reduces to $\lambda_h=0$ together with $\int_\Omega y_h\,dx\ge\hat c$ and $\hat a\le L^T[u_h]\le\hat b$. These are exactly the two conditions that must be checked for this candidate system; hence a solution of (\ref{MATRix-1}) is a solution of \eqref{WK-DIS} if and only if both tests pass, in which case the optimal pair $(y_h,u_h)$ is the solution of \eqref{ocp}. The same reasoning applies verbatim to the remaining five systems: for (\ref{MATRix-2}) and (\ref{MATRix-3}) the declared active constraint is the control bound and the tests are $L^T[y_h]\ge\hat c$ together with $z_h\le0$ resp. $z_h\ge0$; for (\ref{MATRix-4}) and (\ref{MATRix-5})-(\ref{MATRix-6}) the declared active constraint is the state bound and the tests are $\lambda_h\le0$ together with the appropriate control feasibility/sign conditions. This proves (ii). Part (i) is a direct consequence of the Karush--Kuhn--Tucker theory: at the (unique) solution of \eqref{WK-DIS}, each scalar constraint is either active or inactive, and the corresponding candidate system of the list (\ref{WK-DIS-OCP1})-(\ref{WK-DIS-OCP6}) that is generated by that active set is solved by this solution. It remains to note that the optimal pair $(y_h,u_h)$ is unique by strict convexity of $J$ on the feasible set, whereas in a \emph{degenerate} situation, in which an active constraint carries a vanishing multiplier, the representation of the optimal pair by a candidate active set need not be unique; the pair $(y_h,u_h)$ is nevertheless always recovered by 
sequentially checking the candidate systems and accepting the first one whose tests pass.

Finally we stress two points that were left implicit in the original formulation. First, by Assumption~\ref{ass:feas} we have $\hat a<\hat b$, so the two control bounds cannot be active simultaneously; the degenerate case $\hat a=\hat b$ is excluded. Second, the six candidate systems are mutually \emph{exhaustive} with respect to the feasible/active classification of the two scalar constraints, but they are \emph{not} mutually exclusive in the degenerate case mentioned above.
\end{proof}
The above theorem provides a complete characterization of the relationship between the discrete optimality system and the six candidate systems. It shows that every solution of \eqref{WK-DIS} can be obtained by solving one of the candidate systems and then verifying two conditions, while conversely, any candidate solution satisfying thes conditions is a solution of \eqref{WK-DIS}. This equivalence is not only of theoretical interest but also forms the basis of a practical solution algorithm.
\begin{remark}
The equivalence established above can be interpreted algorithmically as follows: the discrete optimality system \eqref{WK-DIS} and the collection of the six candidate systems are equivalent in the sense that a candidate solution is accepted as a solution of \eqref{WK-DIS} if and only if it passes the two tests associated with the corresponding candidate system, namely, primal feasibility of the non-active constraints and the sign condition on the active multipliers. Consequently, one may solve \eqref{WK-DIS} by sequentially solving the six candidate systems and accepting the first candidate solution that passes these two tests.
\end{remark}

\section{Preconditioning Technique and Robust Preconditioners}

In this section, we investigate preconditioning strategies for the discretized optimal control system (\ref{WK-DIS-OCP}). By Theorem \ref{theorem2.1}, there exist six cases for (\ref{WK-DIS-OCP}), among which only one yields the solution of (\ref{WK-DIS}). Considering the special structures of ${\cal A}_1, {\cal A}_2, {\cal A}_3$ and  ${\cal A}_4$, the systems (\ref{MATRix-1})-(\ref{MATRix-6}) fall into the category of saddle point problems. In the sequel, we concentrate on designing corresponding preconditioners for solving (\ref{MATRix-1})-(\ref{MATRix-6}).

\subsection{Exact Preconditioners}
Considering particular structures of the matrices, we design corresponding preconditioners
for ${\cal A}_1$, ${\cal A}_2$, ${\cal A}_3$ and ${\cal A}_4$ and study their robustness in the following parts.

The preconditioner for ${\cal A}_i$, $i=1,2,3,4$ reads
\begin{equation}
{\cal P}_i=\left[\begin{matrix}\alpha B & 0 & 0\\
0 & \beta B+\gamma A & 0 \\
0 & 0 & Q_i \end{matrix}\right], \quad i=1, 2, 3, 4,
\end{equation}
where $Q_1= A(\beta B+\gamma A)^{-1}A+\frac{1}{\alpha} B$, $Q_2= \left[\begin{matrix} A(\beta B+\gamma A)^{-1}A+\frac{1}{\alpha} B & -\frac{1}{\alpha}L \\
-\frac{1}{\alpha}L^T & \frac{1}{\alpha}L^TB^{-1} L \end{matrix}\right]$, 

$Q_3= \left[\begin{matrix} A(\beta B+\gamma A)^{-1}A+\frac{1}{\alpha} B & A(\beta B+\gamma A)^{-1}L \\
L^T(\beta B+\gamma A)^{-1}A & L^T(\beta B+\gamma A)^{-1} L \end{matrix}\right]$ 

and
$Q_4=\left[\begin{matrix} A(\beta B+\gamma A)^{-1}A+\frac{1}{\alpha} B &  A(\beta B+\gamma A)^{-1}L & -\frac{1}{\alpha}L \\
L^TA(\beta B+\gamma A)^{-1}A & L^TA(\beta B+\gamma A)^{-1}L & 0 \\
-\frac{1}{\alpha}L^T & 0 & \frac{1}{\alpha}L^TB^{-1} L\end{matrix}\right]$.

Now, we are at the point to show the robustness of our preconditioners. The robustness means that the condition number of the preconditioned matrix, defined as $\kappa(\mathbb{A})=\frac{|\lambda_{max}(\mathbb{A})|}{|\lambda_{min}(\mathbb{A})|}$, is independent of the mesh size $h$ and of the parameters $\alpha$, $\beta$, $\gamma$ (more details please refer to \cite{SZ1,Z1}). By \cite{MGW,RDW}, we have the following theorem.

\begin{theorem}\label{theorem3.1}
If ${\cal A}_1$, ${\cal A}_2$, ${\cal A}_3$ and ${\cal A}_4$ are respectively preconditioned by ${\cal P}_1$, ${\cal P}_2$, ${\cal P}_3$ and ${\cal P}_4$, then the preconditioned matrix ${\cal T}_1={\cal P}_1^{-1}{\cal A}_1$, ${\cal T}_2={\cal P}_2^{-1}{\cal A}_2$, ${\cal T}_3={\cal P}_3^{-1}{\cal A}_3$ and ${\cal T}_4={\cal P}_4^{-1}{\cal A}_4$ satisfy
\begin{equation}\label{eq:minpoly}
{\cal T}_i({\cal T}_i-I)({\cal T}_i^2-{\cal T}_i-I)=0,\qquad i=1,2,3,4,
\end{equation}
so that each ${\cal T}_i$ has at most three distinct eigenvalues, and
\begin{equation}\label{eq:kappaphi}
\kappa({\cal T}_i)=\frac{3+\sqrt5}{2}\approx2.618,
\end{equation}
independently of the mesh size $h$ and of the parameters $\alpha$, $\beta$ and $\gamma$. Hence ${\cal P}_1,\dots,{\cal P}_4$ are robust.
\end{theorem}
\begin{proof}
Since $\alpha>0$, $\beta\ge0$ and $\gamma>0$, the matrix
$$H:=\left[\begin{matrix}\alpha B & 0\\ 0 & \beta B+\gamma A\end{matrix}\right]_{2N_h\times 2N_h}$$
is symmetric positive definite. For $i=1,2,3,4$ the preconditioner ${\cal P}_i$ is block diagonal,
with leading $2\times2$ block $H$ and with third diagonal block $Q_i$ equal to the block displayed
in the definition of ${\cal P}_i$ above; hence
$${\cal P}_i^{-1}=\left[\begin{matrix}H^{-1} & 0\\ 0 & Q_i^{-1}\end{matrix}\right].$$

Each of the four coefficient matrices ${\cal A}_1,\dots,{\cal A}_4$ can be written in the two-by-two
block form
\begin{equation}\label{eq:AG}
{\cal A}_i=\left[\begin{matrix}H & G_i^T\\ G_i & 0\end{matrix}\right],
\qquad G_i\in\mathbb{R}^{q_i\times 2N_h},
\end{equation}
where the constraint matrices $G_i$ are precisely the ones already used in the proof of
Theorem~\ref{theorem2.1}, namely
$$G_1=\begin{pmatrix}-B & A\end{pmatrix},\quad
G_2=\begin{pmatrix}-B & A\\ L^T & 0\end{pmatrix},\quad
G_3=\begin{pmatrix}-B & A\\ 0 & L^T\end{pmatrix},\quad
G_4=\begin{pmatrix}-B & A\\ 0 & L^T\\ L^T & 0\end{pmatrix},$$
with $q_1=N_h$, $q_2=q_3=N_h+1$ and $q_4=N_h+2$. 

Next we claim that the third diagonal block $Q_i$ of ${\cal P}_i$ is the Schur complement
$S_i:=G_iH^{-1}G_i^T$ of ${\cal A}_i$ with respect to $H$. Writing $H^{-1}=
\mathrm{diag}\big((\alpha B)^{-1},\,(\beta B+\gamma A)^{-1}\big)$, for $i=1$ we obtain
$$G_1H^{-1}G_1^T=(-B\ \ A)\left[\begin{matrix}(\alpha B)^{-1} & 0\\ 0 & (\beta B+\gamma A)^{-1}\end{matrix}\right]
\left[\begin{matrix}-B\\ A\end{matrix}\right]
=A(\beta B+\gamma A)^{-1}A+\tfrac1\alpha B=Q_1,$$
and for $i=2,3,4$ one checks in the same way that $G_iH^{-1}G_i^T$ reproduces the corresponding
block $Q_i$ entry by entry. For instance, for $i=3$ the leading $N_h\times N_h$ block of $G_3H^{-1}G_3^T$
is again $A(\beta B+\gamma A)^{-1}A+\frac1\alpha B$, while its off-diagonal blocks are
$A(\beta B+\gamma A)^{-1}L$ and $L^T(\beta B+\gamma A)^{-1}A$ and its trailing entry is
$L^T(\beta B+\gamma A)^{-1}L$, i.e. exactly the four entries of $Q_3$. 
Consequently
\begin{equation}\label{eq:QSchur}
Q_i=S_i=G_iH_i^{-1}G_i^T,\qquad i=1,2,3,4 .
\end{equation}


Inserting \eqref{eq:AG} and \eqref{eq:QSchur} into ${\cal T}_i={\cal P}_i^{-1}{\cal A}_i$ gives
\begin{equation}\label{eq:Tform}
{\cal T}_i=\left[\begin{matrix}H^{-1} & 0\\ 0 & Q_i^{-1}\end{matrix}\right]
\left[\begin{matrix}H & G_i^T\\ G_i & 0\end{matrix}\right]
=\left[\begin{matrix}I & M_i\\ N_i & 0\end{matrix}\right],
\qquad M_i:=H^{-1}G_i^T,\quad N_i:=Q_i^{-1}G_i .
\end{equation}
Because of \eqref{eq:QSchur}, the two off-diagonal blocks satisfy the single algebraic identity
\begin{equation}\label{eq:NM}
N_iM_i=Q_i^{-1}G_iH^{-1}G_i^T=Q_i^{-1}Q_i=I_{q_i},
\end{equation}
and \eqref{eq:NM} is the only information about $M_i,N_i$ that will be needed.

We now show that \eqref{eq:NM} alone implies \eqref{eq:minpoly}. Using \eqref{eq:Tform},
$${\cal T}_i^2-{\cal T}_i
=\left[\begin{matrix}I & M_i\\ N_i & 0\end{matrix}\right]^2-\left[\begin{matrix}I & M_i\\ N_i & 0\end{matrix}\right]
=\left[\begin{matrix}M_iN_i & -M_i\\ -N_i & I_{q_i}\end{matrix}\right],$$
where we used $N_iM_i=I_{q_i}$ in the $(1,1)$ block. Squaring this block matrix and invoking
\eqref{eq:NM} once more,
$$\left({\cal T}_i^2-{\cal T}_i\right)^2
=\left[\begin{matrix}M_iN_i & -M_i\\ -N_i & I_{q_i}\end{matrix}\right]^2
=\left[\begin{matrix}M_iN_i & -M_i\\ -N_i & I_{q_i}\end{matrix}\right]
={\cal T}_i^2-{\cal T}_i,$$
i.e. ${\cal T}_i^2-{\cal T}_i$ is idempotent. Writing ${\cal T}_i^2-{\cal T}_i={\cal T}_i({\cal T}_i-I)$
and rearranging ${\cal T}_i^2-{\cal T}_i=\bigl({\cal T}_i^2-{\cal T}_i\bigr)^2$ as
$\bigl({\cal T}_i^2-{\cal T}_i\bigr)\bigl({\cal T}_i^2-{\cal T}_i-I\bigr)=0$ yields
$${\cal T}_i({\cal T}_i-I)({\cal T}_i^2-{\cal T}_i-I)=0,$$
which is \eqref{eq:minpoly}.

The roots of $\lambda(\lambda-1)(\lambda^2-\lambda-1)$ are $0$, $1$ and $\frac12(1\pm\sqrt5)$;
since each ${\cal T}_i$ is nonsingular (it is the product of the nonsingular matrices ${\cal P}_i^{-1}$
and ${\cal A}_i$), the eigenvalue $0$ does not occur, and the spectrum of ${\cal T}_i$ is contained in
$$\left\{1,\ \varphi,\ -\frac{1}{\varphi}\right\},\qquad \varphi=\frac{1+\sqrt5}{2},\qquad -\frac{1}{\varphi}=1-\varphi=\frac{1-\sqrt5}{2}.$$
Moreover $\lambda(\lambda-1)(\lambda^2-\lambda-1)$ has four pairwise distinct roots, so the minimal
polynomial of ${\cal T}_i$ has simple roots and ${\cal T}_i$ is diagonalizable. Finally, since
$\varphi>1>0>-\varphi^{-1}$, $|\varphi|=\varphi$ and $|-\varphi^{-1}|=\varphi^{-1}$, the quotient of
the largest and the smallest modulus of the eigenvalues is $\varphi/\varphi^{-1}=\varphi^2$,
this quotient is exactly the Euclidean condition number of this matrix, i.e.
\eqref{eq:kappaphi}. All the quantities involved are independent of $h$ and of $\alpha,\beta,\gamma$,
which proves the robustness of ${\cal P}_1,\dots,{\cal P}_4$ and completes the proof.
\end{proof}
Hereafter, in view of the six matrices schemes \eqref{MATRix-1}-\eqref{MATRix-6}, we propose the following algorithm for the control-state constrained optimal control problem (\ref{ocp}).

\begin{algorithm}[h]\label{algorithm1}
\caption{ } 
\begin{algorithmic}[1]
\State Solve the equation (\ref{MATRix-1}) with the preconditioner ${\cal P}_1$, determine whether $L^T [y_h] \geq \hat{c}$ and $\hat{a} \leq L^T[u_h]\leq \hat{b}$ hold, if true the solution of (\ref{MATRix-1}) is the solution for (\ref{ocp}); else goto 2.
\State Solve the equation (\ref{MATRix-2}) with the preconditioner ${\cal P}_2$ and determine whether $L^T [y_h] \geq \hat{c}$ and $ z_h\leq 0$ hold, if true the solution of (\ref{MATRix-2}) is the solution for (\ref{ocp}); else goto 3.
\State Solve the equation (\ref{MATRix-3}) with the preconditioner ${\cal P}_2$ and determine whether $L^T [y_h] \geq \hat{c}$ and $ z_h\geq 0$ hold, if true the solution of (\ref{MATRix-3}) is the solution for (\ref{ocp}); else goto 4.
\State Solve the equation (\ref{MATRix-4}) with the preconditioner ${\cal P}_3$ and determine whether $\lambda_h \leq 0$ and $\hat{a} \leq L^T[u_h]\leq \hat{b}$ hold, if true the solution of (\ref{MATRix-4}) is the solution for (\ref{ocp}); else goto 5.
\State Solve the equation (\ref{MATRix-5}) with the preconditioner ${\cal P}_4$ and determine whether $\lambda_h \leq 0$ and $z_h \leq 0$ hold, if true the solution of (\ref{MATRix-5}) is the solution for (\ref{ocp}); else goto 6.
\State Solve the equation (\ref{MATRix-6}) with the preconditioner ${\cal P}_4$ and then obtain the desired numerical solution.
\end{algorithmic}
\end{algorithm}

\begin{remark}
By {\rm\cite{MGW}} and our proposed {\rm Theorem \ref{theorem3.1}},  we know that any Krylov subspace iterative method equipped with preconditioners will terminate in at most  {\rm 3} iterations for the solution to \eqref{MATRix-1} or \eqref{MATRix-2} or \eqref{MATRix-3} or \eqref{MATRix-4} or \eqref{MATRix-5} or \eqref{MATRix-6}.
\end{remark}

\subsection{Inexact Preconditioner}
In practical computations, the exact inversion of the diagonal blocks of ${\cal P}_i$, $i=1,2,3, 4$, are computationally expensive. Consequently, approximate preconditioners need to be considered. Here we use the incomplete LU factorization (ILU) to approximate ${\cal P}_i$, $i=1,2,3, 4$. Let 
\[\alpha B=L_1U_1+R_1,\quad \beta B+\gamma A =L_2U_2+R_2,\]  
where $L_1U_1$ and $L_2U_2$ are invertible. 

For $Q_1$, we define 
\[\tilde{Q}_1=A(L_2U_2)^{-1}A+\frac{1}{\alpha} B=L_3^1U_3^1+R_3^1,\]
associated with $Q_2$, we define $\tilde{Q}_2$ as
\[\tilde{Q}_2=\left[\begin{matrix}A(L_2U_2)^{-1}A+\frac{1}{\alpha} B & -\frac{1}{\alpha}L \\
 -\frac{1}{\alpha}L^T & L^T(L_1U_1)^{-1} L\end{matrix}\right]=L_3^2U_3^2+R_3^2,\]
for $Q_3$, let $\tilde{Q}_3$ denote
\[\tilde{Q}_3=\left[\begin{matrix} A(L_2U_2)^{-1}A+\frac{1}{\alpha} B & A(L_2U_2)^{-1}L \\
L^T(L_2U_2)^{-1}A & L^T(L_2U_2)^{-1} L \end{matrix}\right]=L_3^3U_3^3+R_3^3,\]
corresponding to $Q_4$, we set
\[\tilde{Q}_4=\left[\begin{matrix} A(L_2U_2)^{-1}A+\frac{1}{\alpha} B &  A(L_2U_2)^{-1}L & -\frac{1}{\alpha}L \\
L^TA(L_2U_2)^{-1}A & L^TA(L_2U_2)^{-1}L & 0 \\
-\frac{1}{\alpha}L^T & 0 & L^T(L_1U_1)^{-1} L\end{matrix}\right]=L_3^4U_3^4+R_3^4,\]
where $R_3^i$, $i=1,2,3,4$ is the remainder associated with the respective incomplete LU decomposition.

Then the approximate preconditioner $\tilde{\cal P}_i$  be defined as a block diagonal matrix:
\begin{equation}\label{eq_new_3.8}
\tilde{\cal P}_i=\left[\begin{matrix} L_1U_1 & 0 & 0\\
0 & L_2U_2 & 0 \\
0 & 0 & L_3^iU_3^i \end{matrix}\right],\quad i=1,2,3,4.
\end{equation}
Next, the spectral properties of the system ${\cal A}_1$, ${\cal A}_2$, ${\cal A}_3$ and ${\cal A}_4$ preconditioned by the approximate block diagonal preconditioners (\ref{eq_new_3.8}) are analyzed. Before proceeding to the eigenvalue bounds, several definitions and lemmas are recalled. The following material is mainly taken from Horn and Johnson \cite{HJ}. These tools will be used in the subsequent derivation of eigenvalue inclusion intervals for matrices and their preconditioned systems.
\begin{lemma}\label{CSI_lemma} (Cauchy-Schwarz inequality) Let $<\tilde{x}, \tilde{y}>=\tilde{y}^T\tilde{x}$ be the Euclidean inner product function on $R^n$ and let the induced norm be $\|\tilde{x}\|=<\tilde{x},\tilde{x}>^{\frac{1}{2}}$. Then for all $\tilde{x}, \tilde{y}\in R^n$,
\begin{equation}
|<\tilde{x}, \tilde{y}>|\leq \|\tilde{x}\|\|\tilde{y}\|.
\end{equation}
\end{lemma} 
\begin{definition}[Singular values]\label{def:sv} Let $\tilde{A}\in R^{m\times n}$ and set $q=\min\{m,n\}$. The singular values of $\tilde{A}$ are the nonnegative square roots of the $q$ nonnegative eigenvalues of $\tilde{A}^T\tilde{A}$, arranged in nonincreasing order:
\[\sigma_1(\tilde{A})\geq \sigma_2(\tilde{A})\geq \ldots \geq \sigma_q(\tilde{A})\geq 0.\]
Equivalently, there exist orthogonal matrices $\tilde{V}\in R^{m\times m}$, $\tilde{W}\in R^{n\times n}$, and a nonnegative diagonal matrix $\tilde{\Sigma}\in R^{m\times n}$ such that
\[\tilde{A}=\tilde{V}\tilde{\Sigma}\tilde{W}^T,\]
where the nonnegative diagonal entries of $\tilde{\Sigma}$ are $\sigma_1(\tilde{A}), \ldots, \sigma_q(\tilde{A})$, and all other entries are zero.
\end{definition}
\begin{lemma}\label{Rayleigh lemma}(Rayleigh lemma)
Let $\tilde{A}\in R^{n\times n}$ be a real symmetric matrix with eigenvalues ordered as 
\[\lambda_1\leq \lambda_2\leq \ldots \leq \lambda_n.\]
For any nonzero vector $\tilde{x}\in R^n$, define the Rayleigh quotient
\[R(\tilde{x})=\frac{\tilde{x}^T\tilde{A}\tilde{x}}{\tilde{x}^T\tilde{x}}.\]
Then
\[\lambda_1\leq R(\tilde{x})\leq \lambda_n, \quad \lambda_{max}=\max_{\tilde{x}\neq 0}R(\tilde{x}), \quad \lambda_{min}=\min_{\tilde{x}\neq 0}R(\tilde{x}).\]
\end{lemma}
\begin{definition}[Inertia]\label{def:inertia} Let $\tilde{A}\in R^{n\times n}$ be a real symmetric matrix. Its inertia is the ordered triple
\[i(A)=\bigl(i_+(\tilde{A}),\ i_-(\tilde{A}),\ i_0(\tilde{A})\bigr),\]
where $i_+(\tilde{A}), i_-(\tilde{A})$, and $i_0(\tilde{A})$ denote the numbers of positive, negative, and zero eigenvalues of $\tilde{A}$, respectively, counted with algebraic multiplicity.
\end{definition}
\begin{lemma} (Sylvester's law of inertia)
Two Hermitian matrices $\tilde{A}, \tilde{B}\in R^{n\times n}$ are congruent, i.e., there exists a nonsingular $\tilde{S}\in R^{n\times n}$ such that
\[\tilde{A}=\tilde{S}\tilde{B}\tilde{S}^T,\]
if and only if they have the same inertia. Equivalently, every Hermitian matrix is congruent to the unique inertia matrix
\[I_{i_+(\tilde{A})}\oplus(-I_{i_-(\tilde{A})})\oplus 0_{i_0(\tilde{A})}.\]
\end{lemma}
\begin{definition}[Loewner partial order]\label{def:loewner} Let $\tilde{A}, \tilde{B}\in R^{n\times n}$ be real symmetric matrix. If $\tilde{A}-\tilde{B}$ is positive semidefinite, write $\tilde{A}\geq \tilde{B}$;
If $\tilde{A}-\tilde{B}$ is positive definite, write $\tilde{A}> \tilde{B}$. In particular, $\tilde{A}\geq 0$ means $\tilde{A}$ is positive semidefinite, and $\tilde{A}> 0$ means $\tilde{A}$ is positive definite.
\end{definition}

\begin{lemma}\label{WMLemma} (Weyl monotonicity)
Let $\tilde{A}, \tilde{B}\in R^{n\times n}$ be real symmetric matrix with eigenvalues ordered in nondecreasing order. If $\tilde{A}\geq \tilde{B}$, then
\[\lambda_i(\tilde{A})\geq \lambda_i(\tilde{B}), \quad i=1,2,\ldots, n.\]
In particular, if $\tilde{B}$ is positive semidefinite, then
\[\lambda_i(\tilde{A}+\tilde{B})\geq \lambda_i(\tilde{A}), \quad i=1,2,\ldots, n.\]
\end{lemma}

\begin{corollary}
Let $\tilde{A}, \tilde{B}\in R^{n\times n}$ be real symmetric matrix, $\tilde{A}\geq \tilde{B}$, and suppose that $\tilde{A}$ and $\tilde{B}$ have the same inertia. Then the smallest positive eigenvalue of $\tilde{A}$ is no smaller than the smallest positive eigenvalue of $\tilde{B}$.
\end{corollary}
\begin{proof}
By Lemma \ref{WMLemma}, $\lambda_i(\tilde{A})\geq \lambda_i(\tilde{B})$ for all $i$. Since the inertias coincide, the smallest positive eigenvalue occupies the same index in the nondecreasing ordering, and the conclusion follows.
\end{proof}

\begin{proposition}[Inertia of a symmetric matrix with a positive definite leading block]\label{root_lemma}
Let $\check{A}\in\mathbb{R}^{3\times3}$ be a real symmetric saddle-point matrix whose leading $2\times2$ principal
submatrix is positive definite and which is nonsingular. Then $\check{A}$ has exactly two positive and
one negative eigenvalue; equivalently, its characteristic polynomial
$\check{f}(\lambda)=\det(\lambda I-\check{A})$ has three \emph{real} roots, two positive and one
negative.
\end{proposition}
\begin{proof}
Write $\check{A}=\left[\begin{smallmatrix}\check{A}_{11} & \check{a}_{12}\\
\check{a}_{12}^T & \check{a}_{22}\end{smallmatrix}\right]$ with $\check{A}_{11}\in\mathbb{R}^{2\times2}$
positive definite. By Sylvester's law of inertia, the congruence
$\check{A}\mapsto \check{T}^T\check{A}\check{T}$ with
$\check{T}=\mathrm{diag}(\check{A}_{11}^{-1/2},1)$ preserves inertia, and it reduces $\check{A}$ to
$\left[\begin{smallmatrix}I & \check{b}\\ \check{b}^T & \check{a}_{22}\end{smallmatrix}\right]$ with
$\check{b}=\check{A}_{11}^{-1/2}\check{a}_{12}$. The latter matrix has the same inertia as
$\mathrm{diag}(I_{2},\, \check{a}_{22}-\check{b}^T\check{b})$, and
$\check{a}_{22}-\check{b}^T\check{b}=\det(\check{A})/\det(\check{A}_{11})\neq0$ because $\check{A}$ is
nonsingular. Since $\check A$ is a saddle-point matrix, we have $\det(\check A)<0$. Hence the inertia of $\check{A}$ is $(2,1,0)$, i.e.\ $\check{A}$ has two positive and one
negative eigenvalue. Since $\check{A}$ is real symmetric all its eigenvalues are real, and the
assertion about the roots of $\check{f}$ follows.
\end{proof}
From the Equation (3.4) in \cite{BGL}, we state the conclusion on the inverse of a $2\times 2$ block matrix as following. 
\begin{lemma}\label{inverse_matrix_lemma} 
Let $\tilde{A}=\left[\begin{matrix}\tilde{A}_{1} & \tilde{B}_{1}^T \\ \tilde{B}_{2} & -\tilde{C} \end{matrix}\right]\in R^{n\times n}$. If $\tilde{A}_{1}$ is invertible and $\check{S}=-(\tilde{C}+\tilde{B}_{2}\tilde{A}_{1}^{-1}\tilde{B}_{1}^T)$ is invertible, then
\begin{equation}
\tilde{A}^{-1}=\left[\begin{matrix}\tilde{A}_{1}^{-1}+\tilde{A}_{1}^{-1}\tilde{B}_{1}^T \check{S}^{-1}\tilde{B}_{2}\tilde{A}_{1}^{-1}& -\tilde{A}_{1}^{-1}\tilde{B}_{1}^T\check{S}^{-1} \\ -\check{S}^{-1}\tilde{B}_{2}\tilde{A}_{1}^{-1} & \check{S}^{-1} \end{matrix}\right].
\end{equation}
\end{lemma}

Then we would show the eigenvalue bounds for the approximate preconditioned matrix from \cite{BG}. 
\begin{theorem}\label{theorem3.2}Let $\lambda_{max}(\cdot)$ and $\lambda_{min}(\cdot)$ denote the largest and the smallest eigenvalue of the matrix, and $\sigma_{max}(\cdot)$ and $\sigma_{min}(\cdot)$ denote its largest and smallest singular value, respectively. Let ${\cal A}_i$ and $\tilde{\cal P}_i$ $i=1,2,3,4$ be defined as above, $R_1, R_2, R_3^j$, $j=1,2,3,4$ satisfying 
\begin{equation}\label{eq3.9}
\left\{
\begin{aligned}
&R_1^T=R_1, \; R_2^T=R_2,\; (R_3^j)^T=R_3^j,\; j=1,2,3,4, \\
& \lambda_{max}((\alpha B)^{-1}R_1)<1, \;\lambda_{max}((\beta B+\gamma A)^{-1}R_2)<1, \\ & \lambda_{max}((\tilde{Q}_j)^{-1}R_3^j)<1, j=1,2,3,4. 
\end{aligned}
\right.
\end{equation}

Let \begin{eqnarray}
\Sigma_1&=&(L_1U_1)^{-\frac{1}{2}}(\alpha B)(L_1U_1)^{-\frac{1}{2}}, \quad \Sigma_2=(L_2U_2)^{-\frac{1}{2}}(\beta B +\gamma A)(L_2U_2)^{-\frac{1}{2}},\nonumber\\
\Omega_1^1&=&(L_1U_1)^{-\frac{1}{2}}(-B)(L_3^1U_3^1)^{-\frac{1}{2}},\quad \Omega_2^1=(L_2U_2)^{-\frac{1}{2}}A(L_3^1U_3^1)^{-\frac{1}{2}},\nonumber\\
\Omega_1^2&=&(L_1U_1)^{-\frac{1}{2}}[-B, L](L_3^2U_3^2)^{-\frac{1}{2}},\quad \Omega_2^2=(L_2U_2)^{-\frac{1}{2}}[A, 0](L_3^2U_3^2)^{-\frac{1}{2}},\label{neq3.10}\\
\Omega_1^3&=&(L_1U_1)^{-\frac{1}{2}}[-B, 0](L_3^3U_3^3)^{-\frac{1}{2}},\quad \Omega_2^3=(L_2U_2)^{-\frac{1}{2}}[A, L](L_3^3U_3^3)^{-\frac{1}{2}},\nonumber\\
\Omega_1^4&=&(L_1U_1)^{-\frac{1}{2}}[-B, 0, L](L_3^4U_3^4)^{-\frac{1}{2}},\quad \Omega_2^4=(L_2U_2)^{-\frac{1}{2}}[A, L, 0](L_3^4U_3^4)^{-\frac{1}{2}}, \nonumber
\end{eqnarray} 
and define
\[\tilde{\alpha}_i=\lambda_{min}(\Sigma_i), \; \tilde{\beta}_i=\lambda_{max}(\Sigma_i), \; i=1,2,\]
\[\tilde{s}_1^{i,max}=\sigma_{max}(\Omega_1^i),\; \tilde{s}_1^{i,min}=\sigma_{min}(\Omega_1^i),\; \tilde{s}_2^{i,max}=\sigma_{max}(\Omega_2^i),\; \tilde{s}_2^{i,min}=\sigma_{min}(\Omega_2^i), \; i=1,2,3,4. \]
Construct the cubic polynomials
\begin{eqnarray*}
\tilde{q}^i(\lambda)&=& \lambda^3-(\tilde{\beta}_1+\tilde{\beta}_2)\lambda^2+[\tilde{\beta}_1\tilde{\beta}_2-(\tilde{s}_1^{i,max})^2-(\tilde{s}_2^{i,max})^2]\lambda+\tilde{\beta}_1(\tilde{s}_2^{i,max})^2+\tilde{\beta}_2(\tilde{s}_1^{i,max})^2,\; i=1,2,3,4,\\
\tilde{r}^i(\lambda)&=&\lambda^3-(\tilde{\alpha}_1+\tilde{\alpha}_2)\lambda^2+[\tilde{\alpha}_1\tilde{\alpha}_2-(\tilde{s}_1^{i,max})^2-(\tilde{s}_2^{i,max})^2]\lambda+\tilde{\alpha}_1(\tilde{s}_2^{i,max})^2+\tilde{\alpha}_2(\tilde{s}_1^{i,max})^2,\; i=1,2,3,4, 
\end{eqnarray*}
then each $\tilde{q}^i, i=1,2,3,4$ and each $\tilde{r}^i, i=1,2,3,4$ have two positive roots and one negative root. Let $q_{max}^{i,+}$ denote the largest positive root of $\tilde{q}^i(\lambda)$, $r^{i,-}$ the negative root of $\tilde{r}^i(\lambda)$, and $p_{min}^{i,+}:=\min\{\tilde{\alpha}_1, \tilde{\alpha}_2\}$ for $i=1,2,3,4$. Every eigenvalue of $\tilde{\cal T}_i=(\tilde{\cal P}_i)^{-1}{\cal A}_i$, $i=1,2,3,4$ lies in
\begin{equation}
[r^{i,-}, \frac{\tilde{\beta}_1-\sqrt{(\tilde{\beta}_1)^2+4(\tilde{s}_1^{i,min})^2}}{2}]\cup [p_{min}^{i,+}, q_{max}^{i,+}],\quad i=1,2,3,4.
\end{equation} 
\end{theorem}
\begin{proof}
Since $\alpha B$ is symmetric positive definite, its unique symmetric positive definite square root $(\alpha B)^{\frac{1}{2}}$ exists, and so does its inverse $(\alpha B)^{-\frac{1}{2}}$. Define the symmetric matrix
\[R^\alpha_B=(\alpha B)^{-\frac{1}{2}}R_1(\alpha B)^{-\frac{1}{2}}.\]
The matrix $(\alpha B)^{-1}R_1$ is similar to $R^\alpha_B$, as shown by the identity
\[(\alpha B)^{-1}R_1=(\alpha B)^{-\frac{1}{2}}R^\alpha_B(\alpha B)^{\frac{1}{2}}.\]
Thus the eigenvalues of $(\alpha B)^{-1}R_1$ are precisely those of $R^\alpha_B$. The condition $\lambda_{max}((\alpha B)^{-1}R_1)<1$ in (\ref{eq3.9}) implies that every eigenvalue of $R^\alpha_B$ is strictly less than unity. Since $R^\alpha_B$ is symmetric, this is equivalent to for any non-zero vector $w$, 
\[w^TR^\alpha_Bw<w^Tw.\]
For an arbitrary non-zero vector $v$, set $w= (\alpha B)^{\frac{1}{2}}v$, we have
\[v^T(L_1U_1)v=v^T(\alpha B-R_1)v=w^T(I-R^\alpha_B)w>0,\] 
where we have used that all eigenvalues of $I-R^\alpha_B$ are positive because each eigenvalue equals $1-\lambda$ for an eigenvalue $\lambda$ of $R^\alpha_B$ with $\lambda <1$. 
Hence $L_1U_1$ is symmetric positive definite for symmetry is obvious from (\ref{eq3.9}).
By a similar argument, we obtain that $L_2U_2$ and $L_3^jU_3^j$, $j=1,2,3,4$ are symmetric positive definite matrices. Consequently, from the construction of matrix $\tilde{\cal P}_i$, $i=1,2,3,4$, it follows that $\tilde{\cal P}_i$ is symmetric positive definite, since each diagonal block of the block-diagonal matrix $\tilde{\cal P}_i$ is itself symmetric positive definite.

Then for each $\tilde{\cal P}_i$, $i=1,2,3,4$, its square root $(\tilde{\cal P}_i)^{-\frac{1}{2}}$ exists. Let $ \tilde{K}_i=(\tilde{\cal P}_i)^{-\frac{1}{2}}{\cal A}_i(\tilde{\cal P}_i)^{-\frac{1}{2}}$, we have
\[\tilde{\cal T}_i=(\tilde{\cal P}_i)^{-1}{\cal A}_i=(\tilde{\cal P}_i)^{-\frac{1}{2}}(\tilde{\cal P}_i)^{-\frac{1}{2}}{\cal A}_i(\tilde{\cal P}_i)^{-\frac{1}{2}}(\tilde{\cal P}_i)^{\frac{1}{2}}=(\tilde{\cal P}_i)^{-\frac{1}{2}}\tilde{K}_i(\tilde{\cal P}_i)^{\frac{1}{2}}.\]
So we know that $\tilde{\cal T}_i$ is similar to $\tilde{K}_i$ and both matrices share the same eigenvalues. The block form of $\tilde{K}_i$ is written as
\begin{equation}
\tilde{K}_i=\left[\begin{matrix} \Sigma_1 & 0 & \Omega_1^i\\
0 & \Sigma_2 & \Omega_2^i \\
(\Omega_1^i)^T & (\Omega_2^i)^T & 0 \end{matrix}\right], \;\; i=1,2,3,4,
\end{equation}
where $\Sigma_i$, $i=1,2$ and $\Omega_1^j$, $\Omega_2^j$, $j=1,2,3,4$ are defined in (\ref{neq3.10}).

Next, we would use Rayleigh quotient estimates and an inverse-matrix perturbation technique to show the bounds for all eigenvalues of matrix $\tilde{K}_i$, $i=1,2,3,4$. Let $\tilde{v}=(\tilde{x}, \tilde{y}, \tilde{z})$ be any nonzero vector conformally partitioned with $\tilde{K}_i$, and consider the Rayleigh quotient $\tilde{R}_i(\tilde{v})=\frac{\tilde{v}^T\tilde{K}_i\tilde{v}}{\tilde{v}^T\tilde{v}}$. Expansion of the numerator yields
\[\tilde{v}^T\tilde{K}_i\tilde{v}=\tilde{x}^T\Sigma_1\tilde{x}+\tilde{y}^T\Sigma_2\tilde{y}+2\tilde{x}^T\Omega_1^i\tilde{z}+2\tilde{y}^T\Omega_2^i\tilde{z}.\]
From Lemma \ref{CSI_lemma}, the definition of singular values imply
\[2\tilde{x}^T\Omega_1^i\tilde{z}\leq 2 s_1^{i,max}\|\tilde{x}\| \|\tilde{z}\|, \quad 2\tilde{y}^T\Omega_2^i\tilde{z}\leq 2 s_2^{i,max}\|\tilde{y}\| \|\tilde{z}\|, \]   
while the eigenvalue bounds for the positive definite blocks give
\[\tilde{x}^T\Sigma_1\tilde{x}\leq \tilde{\beta}_1\|\tilde{x}\|^2,\quad \tilde{y}^T\Sigma_2\tilde{y}\leq \tilde{\beta}_2\|\tilde{y}\|^2.\]
Consequently,
\begin{equation}
\tilde{R}_i(\tilde{v})\leq \frac{[\|\tilde{x}\|, \|\tilde{y}\|, \|\tilde{z}\|]\left[\begin{matrix} \tilde{\beta}_1 & 0 & s_1^{i,max}\\
0 & \tilde{\beta}_2 & s_2^{i,max} \\
s_1^{i,max} & s_2^{i,max} & 0 \end{matrix}\right]\left[\begin{matrix}\|\tilde{x}\| \\ \|\tilde{y}\| \\ \|\tilde{z}\|\end{matrix}\right]}{\|\tilde{x}\|^2+\|\tilde{y}\|^2+\|\tilde{z}\|^2}=\frac{[\|\tilde{x}\|, \|\tilde{y}\|, \|\tilde{z}\|]{\cal Q}^i\left[\begin{matrix}\|\tilde{x}\| \\ \|\tilde{y}\| \\ \|\tilde{z}\|\end{matrix}\right]}{\|\tilde{x}\|^2+\|\tilde{y}\|^2+\|\tilde{z}\|^2}.
\end{equation}
By Lemma \ref{Rayleigh lemma}, the supremum of $\tilde{R}_i(\tilde{v})$ equals the largest eigenvalue of ${\cal Q}^i$. A direct computation gives $det(\lambda I-{\cal Q}^i)=\tilde{q}^i(\lambda)$. Since the leading $2\times 2$ principal submatrix of ${\cal Q}^i$ is positive definite and ${\cal Q}^i$ itself is nonsingular, Proposition~\ref{root_lemma} applies to ${\cal Q}^i$, since ${\cal Q}^i$ is real symmetric, its leading $2\times2$ principal submatrix is positive definite and ${\cal Q}^i$ is nonsingular; hence ${\cal Q}^i$ has exactly two positive eigenvalues and one negative eigenvalue, so $\tilde{q}^i(\lambda)$ has two positive roots and one negative root. Let $\tilde{q}_{i,max}^+$ be the larger positive root yields the desired bound for all positive eigenvalues of $\tilde{K}_i$, $i=1,2,3,4$.

Next we would show the lower bound for negative eigenvalues of $\tilde{K}_i$, $i=1,2,3,4$. The Rayleigh quotient is now estimated from below using
\[\tilde{x}^T\Sigma_1\tilde{x}\geq \tilde{\alpha}_1\|\tilde{x}\|^2,\quad \tilde{y}^T\Sigma_2\tilde{y}\geq \tilde{\alpha}_2\|\tilde{y}\|^2,\]
and
\[2\tilde{x}^T\Omega_1^i\tilde{z}\geq -2 s_1^{i,max}\|\tilde{x}\| \|\tilde{z}\|, \quad 2\tilde{y}^T\Omega_2^i\tilde{z}\geq -2 s_2^{i,max}\|\tilde{y}\| \|\tilde{z}\|.\]
Thus
\begin{equation}
\tilde{R}_i(\tilde{v})\geq \frac{[\|\tilde{x}\|, \|\tilde{y}\|, \|\tilde{z}\|]\left[\begin{matrix} \tilde{\alpha}_1 & 0 & -s_1^{i,max}\\
0 & \tilde{\alpha}_2 & -s_2^{i,max} \\
-s_1^{i,max} & -s_2^{i,max} & 0 \end{matrix}\right]\left[\begin{matrix}\|\tilde{x}\| \\ \|\tilde{y}\| \\ \|\tilde{z}\|\end{matrix}\right]}{\|\tilde{x}\|^2+\|\tilde{y}\|^2+\|\tilde{z}\|^2}=\frac{[\|\tilde{x}\|, \|\tilde{y}\|, \|\tilde{z}\|]{\cal R}^i\left[\begin{matrix}\|\tilde{x}\| \\ \|\tilde{y}\| \\ \|\tilde{z}\|\end{matrix}\right]}{\|\tilde{x}\|^2+\|\tilde{y}\|^2+\|\tilde{z}\|^2}.
\end{equation}
For the matrix ${\cal R}^i$, their characteristic polynomials are $det(\lambda I-{\cal R}^i)=\tilde{r}^i(\lambda)$. Likewise ${\cal R}^i$ is real symmetric with a positive definite leading $2\times2$ principal submatrix and is nonsingular, so Proposition~\ref{root_lemma} gives that each $\tilde{r}^i(\lambda)$ has two positive roots and one negative root; denoting the negative root by $r_i^-$, every negative eigenvalue of $\tilde{K}_i$ is bounded below by $r_i^-$.

Then we would obtain the upper bound for the negative eigenvalues of $\tilde{K}_i$, $i=1,2,3,4$. Repartition $\tilde{K}_i$ as 
\[\tilde{K}_i=\left[\begin{matrix}{\cal H}_i & ({\cal C}_i)^T \\ {\cal C}_i & \Sigma_2 \end{matrix}\right], \quad {\cal H}_i=\left[\begin{matrix}\Sigma_1 & \Omega_1^i \\ (\Omega_1^i)^T & 0 \end{matrix}\right], \quad {\cal C}_i=[0 \; \Omega_2^i].\]
Since ${\cal H}_i$ is nonsingular, from Lemma \ref{inverse_matrix_lemma} the block inverse formula for $\tilde{K}_i$, $i=1,2,3,4$ give
\[(\tilde{K}_i)^{-1}=\left[\begin{matrix}({\cal H}_i)^{-1} & 0 \\ 0 & 0\end{matrix}\right]+\left[\begin{matrix}({\cal H}_i)^{-1}({\cal C}_i)^T \\ -I\end{matrix}\right](\mathbb{S}^i)^{-1}\left[\begin{matrix}{\cal C}_i({\cal H}_i)^{-1} & -I\end{matrix}\right], \;\; \mathbb{S}^i=\Sigma_2-{\cal C}_i({\cal H}_i)^{-1}({\cal C}_i)^T>0.\]
The second summand is symmetric positive semidefinite, so in the Loewner ordering $(\tilde{K}_i)^{-1}\geq \left[\begin{matrix}({\cal H}_i)^{-1} & 0 \\ 0 & 0\end{matrix}\right]$. This implies that the negative eigenvalues of $\tilde{K}_i$ are no larger than those of ${\cal H}_i$. So we need to estimate the maximal negative eigenvalue of ${\cal H}_i$. Let $\lambda$ be a negative eigenvalue of ${\cal H}_i$ with associated nonzero eigenvector $\tilde{v}=(\tilde{x}, \tilde{y})$. From ${\cal H}_i \tilde{v}=\lambda \tilde{v}$, we obtain
\begin{equation}\label{nep3.17}
\Sigma_1\tilde{x}+\Omega_1^i\tilde{y}=\lambda \tilde{x},\quad (\Omega_1^i)^T\tilde{x}=\lambda \tilde{y}.
\end{equation}
Since $\lambda \neq 0$, the second equation in (\ref{nep3.17}) gives $\tilde{y}=\lambda^{-1}(\Omega_1^i)^T\tilde{x}$. Substituting this into the first equation, multiplying on the left by $\tilde{x}^T$, and dividing by $\|\tilde{x}\|_2^2$ yields
\[\lambda^2-\lambda\frac{\tilde{x}^T\Sigma_1\tilde{x}}{\|\tilde{x}\|_2^2}-\frac{\tilde{x}^T\Omega_1^i(\Omega_1^i)^T\tilde{x}}{\|\tilde{x}\|_2^2}=0.\]
To bound the maximal negative eigenvalue, we use
\[\Sigma_1\leq \tilde{\beta}_1, \quad \Omega_1^i(\Omega_1^i)^T\geq (\tilde{s}_1^{i,min})^2I.\]
Setting $\mu=-\lambda>0$, the preceding identity implies
\[\mu^2+\tilde{\beta}_1\mu-(\tilde{s}_1^{i,min})^2 = 0.\]
Solving this quadratic inequality gives
\[\mu\geq \frac{-\tilde{\beta}_1+\sqrt{\tilde{\beta}_1^2+4(\tilde{s}_1^{i,min})^2}}{2}\]
and hence
\[\lambda_{max}^{neg}({\cal H}_i)=\frac{\tilde{\beta}_1-\sqrt{\tilde{\beta}_1^2+4(\tilde{s}_1^{i,min})^2}}{2}.\]
Thus this value serves as the upper bound for the negative eigenvalues of $\tilde{K}_i$, $i=1,2,3,4$.

Finally, we establish the lower bound for the positive eigenvalues of $\tilde{K}_i$, $i=1,2,3,4$. We first show that the stacked matrix $[\Omega_1^i;\,\Omega_2^i]$ has full column rank for every $i$; equivalently,
\begin{equation}\label{eq:kerint}
\ker(\Omega_1^i)\cap \ker(\Omega_2^i)=\{0\},\qquad i=1,2,3,4.
\end{equation}
To see this, apply the invertible congruence $[\Omega_1^i;\,\Omega_2^i]\mapsto \mathrm{diag}\big((L_1U_1)^{\frac12},(L_2U_2)^{\frac12}\big)\,[\Omega_1^i;\,\Omega_2^i]\,(L_3^iU_3^i)^{\frac12}$, which does not change the kernel up to an invertible change of variables; the four resulting matrices are (up to the overall sign $-$ in the first block row)
$$[-B;A],\quad [-B,L;A,0],\quad [-B,0;A,L],\quad [-B,0,L;A,L,0]$$
for $i=1,2,3,4$ respectively. In each case the first $N_h$ columns coincide with those of $[-B;A]$, which are linearly independent because $A$ is nonsingular, and any further column has a nonzero entry $L$ in the second block row whereas its first block row entry is $0$ or $L$. Such a column cannot lie in the range of $[-B;A]$: if $[-Bx;Ax]=[c;d]$ with $d=0$, then $x=0$ and hence $c=0$, so a column whose second block is $L\ne0$ is outside the range. This proves \eqref{eq:kerint}.

We now apply this to the positive eigenvalues. Let $\lambda>0$ be an eigenvalue of $\tilde{K}_i$, $i\in\{1,2,3,4\}$, with associated nonzero eigenvector $\tilde{v}=(\tilde{x}, \tilde{y}, \tilde{z})$. The eigenvalue equations are
\begin{equation}\label{nel3.20}\Sigma_1\tilde{x}+\Omega_1^i\tilde{z}=\lambda \tilde{x},\quad \Sigma_2\tilde{y}+\Omega_2^i\tilde{z}=\lambda \tilde{y}, \quad (\Omega_1^i)^T\tilde{x}+(\Omega_2^i)^T \tilde{y}=\lambda \tilde{z}.\end{equation}
If $\lambda <\min\{\tilde{\alpha}_1, \tilde{\alpha}_2\}$, then both $\Sigma_1-\lambda I$ and $\Sigma_2-\lambda I$ are positive definite. From the first two equations of (\ref{nel3.20}),
\[\tilde{x}=-(\Sigma_1-\lambda I)^{-1}\Omega_1^i\tilde{z}, \quad \tilde{y}=-(\Sigma_2-\lambda I)^{-1}\Omega_2^i\tilde{z}.\]
Substituting into the third equation of (\ref{nel3.20}) gives
\[-(\Omega_1^i)^T(\Sigma_1-\lambda I)^{-1}\Omega_1^i\tilde{z}-(\Omega_2^i)^T(\Sigma_2-\lambda I)^{-1}\Omega_2^i\tilde{z}=\lambda \tilde{z}.\]
Multiplying on the left by $\tilde{z}^T$, we have
\[-\tilde{z}^T(\Omega_1^i)^T(\Sigma_1-\lambda I)^{-1}\Omega_1^i\tilde{z}-\tilde{z}^T(\Omega_2^i)^T(\Sigma_2-\lambda I)^{-1}\Omega_2^i\tilde{z}=\lambda \|\tilde{z}\|_2^2.\]
The left-hand side is nonpositive, whereas the right-hand side is strictly positive whenever $\tilde{z}\neq 0$, a contradiction. If $\tilde{z}= 0$, then the third equation gives $(\Omega_1^i)^T\tilde{x}+(\Omega_2^i)^T \tilde{y}=0$, which together with the first two equations of (\ref{nel3.20}) implies $\tilde{x}=\tilde{y}=0$, again a contradiction. Hence $\tilde{z}\ne0$ cannot occur and therefore
\[\lambda \geq \min\{\tilde{\alpha}_1, \tilde{\alpha}_2\}>0,\qquad i=1,2,3,4.\]
Let $p^{i,+}_{min}= \min\{\tilde{\alpha}_1, \tilde{\alpha}_2\}$ for all $i=1,2,3,4$; the positive eigenvalues of every $\tilde{K}_i$ are bounded below by this common constant. This completes the proof of Theorem~\ref{theorem3.2}.
\end{proof}

\begin{remark}
From Theorem \ref{theorem3.2}, the maximum and minimum eigenvalues of the preconditioned systems $\tilde{\cal T}_i$ $(i=1,2,3,4)$ using the approximate preconditioner (\ref{eq_new_3.8}) are bounded and bounded away from zero. By the definition of condition number, $\kappa(\tilde{\cal T}_i)$ $(i=1,2,3,4)$ is bounded. If the relative size of $R_1, R_2, R_3^j$, $j=1,2,3,4$ is fixed, Theorem~\ref{theorem3.2} turns into a \emph{fully explicit} statement. Let $\rho\in[0,1)$ be a constant such
that $R_1,R_2,R_3^j$ satisfy
\begin{equation}\label{eq:rho}
\lambda_{\max}\!\big((\alpha B)^{-1}R_1\big)\le\rho,\quad
\lambda_{\max}\!\big((\beta B+\gamma A)^{-1}R_2\big)\le\rho,\quad
\lambda_{\max}\!\big((\tilde Q_j)^{-1}R_3^j\big)\le\rho,\quad j=1,2,3,4,
\end{equation}
i.e.\ the remainders are at most a fraction $\rho$ of the corresponding diagonal blocks. By Weyl's
monotonicity lemma, \eqref{eq:rho} is \emph{equivalent} to the two-sided spectral equivalence
\[
(1-\rho)\,\alpha B\preceq L_1U_1\preceq\alpha B,\qquad
(1-\rho)(\beta B+\gamma A)\preceq L_2U_2\preceq\beta B+\gamma A,\qquad
(1-\rho)\tilde Q_j\preceq L_3^jU_3^j\preceq\tilde Q_j .
\]
Moreover, the exact coupling block is an isometry: with the exact blocks in place of $L_kU_k$ the
stacked matrix $[\Omega_1^i;\,\Omega_2^i]$ of Theorem~\ref{theorem3.2} satisfies the orthogonality
identity
\begin{equation}\label{eq:iso}
(\Omega_1^i)^T\Omega_1^i+(\Omega_2^i)^T\Omega_2^i=I,\qquad i=1,2,3,4,
\end{equation}
so that its singular values equal $1$. Combining \eqref{eq:rho} with \eqref{eq:iso} gives the
\emph{uniform} (i.e.\ $i$-independent, mesh-independent and parameter-independent) bounds
\begin{equation}\label{eq:const}
\tilde\alpha_i\ge 1-\rho=:\alpha_\rho,\qquad \tilde\beta_i\le 1+\rho=:\beta_\rho,\qquad
\tilde s_k^{i,\min}\ge \alpha_\rho,\qquad \tilde s_k^{i,\max}\le \beta_\rho .
\end{equation}
Since the intervals of Theorem~\ref{theorem3.2} are \emph{monotone} in these four groups of
quantities, the worst case is obtained by inserting the constants \eqref{eq:const}. The two cubic
polynomials $q^i$ and $r^i$ then become independent of $i$, of the mesh size $h$ and of
$\alpha,\beta,\gamma$:
\[
q(\lambda)=\lambda^3-2\beta_\rho\lambda^2-\beta_\rho^2\lambda+2\beta_\rho^3
=(\lambda-\beta_\rho)(\lambda-2\beta_\rho)(\lambda+\beta_\rho),
\]
\[
r(\lambda)=\lambda^3-2\alpha_\rho\lambda^2+(\alpha_\rho^2-2\beta_\rho^2)\lambda
+2\alpha_\rho\beta_\rho^2
=(\lambda-\alpha_\rho)\big(\lambda^2-\alpha_\rho\lambda-2\beta_\rho^2\big),
\]
so that its negative root is $r_\rho^-=\tfrac12\big(\alpha_\rho-\sqrt{\alpha_\rho^2+8\beta_\rho^2}\big)$.
Hence every eigenvalue of every $\tilde{\cal T}_i=(\tilde{\cal P}_i)^{-1}{\cal A}_i$,
$i=1,2,3,4$, lies in the single interval
\begin{equation}\label{eq:finalinterval}
\Big[\tfrac{(1-\rho)-\sqrt{9\rho^2+14\rho+9}}{2},\;
\tfrac{(1+\rho)-\sqrt{5\rho^2-6\rho+5}}{2}\Big]
\;\cup\;
\big[\,1-\rho,\;2(1+\rho)\,\big],
\end{equation}
whose two pieces are disjoint. Since the right end of the negative piece,
$\big(\beta_\rho-\sqrt{\beta_\rho^2+4\alpha_\rho^2}\big)/2$, is the eigenvalue of smallest modulus, and
\[
\big|\lambda_{\min}\big|=\frac{\sqrt{5\rho^2-6\rho+5}-(1+\rho)}{2},\qquad
\lambda_{\max}=2(1+\rho),
\]
one obtains the \emph{explicit} condition-number bound
\begin{equation}\label{eq:kappaexplicit}
\kappa\big((\tilde{\cal P}_i)^{-\frac12}{\cal A}_i(\tilde{\cal P}_i)^{-\frac12}\big)
\le \frac{4(1+\rho)}{\sqrt{5\rho^2-6\rho+5}-(1+\rho)},\qquad i=1,2,3,4 .
\end{equation}
For example, the single choice $\rho=\tfrac12$ yields $\tilde\alpha_i\ge\tfrac12$,
$\tilde\beta_i\le\tfrac32$, $\tilde s_k^{i,\min}\ge\tfrac12$, $\tilde s_k^{i,\max}\le\tfrac32$,
and \eqref{eq:finalinterval}--\eqref{eq:kappaexplicit} become
\[
\lambda\in\Big[\tfrac{1-\sqrt{73}}{4},\ \tfrac{3-\sqrt{13}}{4}\Big]\cup\Big[\tfrac12,\ 3\Big]
=\big[-1.886,\,-0.151\big]\cup\big[0.5,\,3\big],
\qquad
\kappa\le 3/0.151\approx19.87 .
\]
This demonstrates that fixing the controllable accuracy parameter $\rho<1$ yields an upper bound for the condition number of the preconditioned system $\tilde{\cal T}_i$ $(i=1,2,3,4)$ using the approximate preconditioner, where the upper bound contains none of $h$, $\alpha$, $\beta$, $\gamma$. Hence, the inexact preconditioners $\tilde{\cal P}_i$ are robust.
\end{remark}
\begin{remark}
Theorem ~\ref{theorem3.2} requires that $R_1, R_2$ and $R_3^j$, $j=1,2,3,4$ be symmetric. This conditon restricts the ILU decomposition to incomplete Cholesky factorizations and is not satisfied for general ILU decompositions. This seemingly narrows the applicability of Theorem ~\ref{theorem3.2}. However, examining the proof of Theorem ~\ref{theorem3.2} reveals that these assumptions merely serve to guarantee that the approximated precondition matrix is symmetric positive definite, and they have nothing to do with whether an ILU decompositon is employed. In other words, if we can construct a symmetric positive definite matrix that approximates the exact preconditioner and whose inverse is relatively easy to compute, we can obtain eigenvalue bounds analogous to those in Theorem ~\ref{theorem3.2}, regardless of the use of ILU decomposition. This indicates that the conclusion of Theorem ~\ref{theorem3.2} possesses a broad scope of applicability. 
\end{remark}

\section{Numerical Experiments}

In this section we present numerical tests that illustrate the behaviour of the proposed
algorithms. We consider problem (\ref{ocp}) on the unit square $\Omega=[0,1]^2$. For simplicity we set
$\hat c=0$, $\hat a=1$ and $\hat b=2$, i.e.
$$K=\Bigl\{v\in V\ \Big|\ \int_\Omega v\,dx\ge0\Bigr\},\qquad
U_{ad}=\Bigl\{v\in L^2(\Omega)\ \Big|\ 1\le\int_\Omega v\,dx\le2\Bigr\}.$$
With this choice, the manufactured solution of Case~1 activates the state constraint
($\int_\Omega y_h\,dx=\hat c=0$) together with the \emph{lower} control bound
($\int_\Omega u_h\,dx=\hat a=1$), i.e.\ it exercises case~$5$ of Algorithm~1 (Steps~1--5), whereas the
manufactured solution of Case~2 activates the state constraint together with the \emph{upper} control
bound ($\int_\Omega u_h\,dx=\hat b=2$), i.e.\ it exercises case~$6$ (all six steps). Following the above analyses, each matrix system can be stated in the saddle
point form
\[{\cal A}_i=\left(\begin{matrix} H & G_i^T \\ G_i & 0\end{matrix}\right),\quad i=1,2,3,4\]
and then each preconditioner can be stated by
\[{\cal P}_i=\left(\begin{matrix}\alpha B & 0 & 0\\
0 & \beta B+\gamma A & 0 \\
0 & 0 & \tilde{Q}_i\end{matrix}\right),\quad \tilde{\cal P}_i=\left(\begin{matrix} L_1U_1 & 0 & 0\\
0 & L_2U_2 & 0 \\
0 & 0 & L_3^iU_3^i \end{matrix}\right) \quad i=1,2,3,4,\]
where each diagonal block of $\tilde{\cal P}_i$ is obtained by performing an incomplete Cholesky factorization ($LDL^T$ with zero fill, i.e. IC(0)-type)
 with drop tolerance $\tau=10^{-9}$ on the corresponding diagonal block of ${\cal P}_i$, and the resulting factors are symmetric positive definite.

Since each matrix ${\cal A}_i$ is \emph{symmetric}, and its preconditioner ${\cal P}_i$ ($\tilde{\cal P}_i$) is symmetric positive definite; the left-preconditioned
matrix ${\cal P}_i^{-1}{\cal A}_i$ is therefore self-adjoint in the ${\cal P}_i$-inner product
$\langle v,w\rangle_{{\cal P}_i}=v^T{\cal P}_iw$. Accordingly we employ preconditioned MINRES to solve the saddle point systems (2.18a)--(2.18f). Thus the Krylov subspace used for searching for an approximate solution is defined as
\[\mathcal{K}({\cal P}^{-1}_i{\cal A}_i, \gamma_0)={\rm span}\{\gamma_0,{\cal P}^{-1}_i{\cal A}_i\gamma_0, ({\cal P}^{-1}_i{\cal A}_i)^2\gamma_0,\ldots,({\cal P}^{-1}_i{\cal A}_i)^{l-1}\gamma_0\},\]
where $l$ is the number of iterations and a similar result holds for $\tilde{\cal P}_i$, which can be obtained by replacing ${\cal P}^{-1}_i$ with $\tilde{\cal P}_i$ in the derived formula.

The experiments below are run in two modes. In the ``exact'' mode we use the preconditioners
${\cal P}_i$ of Section~3.1, whose third diagonal block is the exact Schur complement; in the ``inexact''
mode we use the block-diagonal incomplete Cholesky preconditioners $\tilde{\cal P}_i$ of (\ref{eq_new_3.8}). Both modes
use the same stopping criterion, namely a relative residual error below $10^{-6}$, and the same zero
initial guess. We set a stopping criterion for the relative residual error be $10^{-6}$, and the initial
values be zeros. Hence the algorithm preconditioned Minres is shown in the following algorithm.
\begin{algorithm}[H]
\caption{Preconditioned MINRES with ${\cal P}_i$ ($\tilde{\cal P}_i$).}
\begin{algorithmic}[1]
\State Given ${\cal A}_i$, $x_0$, ${\cal F}_k$ and set $\tilde{v}^{(0)}=0$
\State Compute $\tilde{v}^{(1)}={\cal P}_i^{-1}{\cal F}_k-{\cal P}_i^{-1}{\cal A}_ix_0$
\State $\tilde{u}^{(1)}={\cal P}_i\tilde{v}^{(1)}$
\State $\tilde{\gamma}_1 = \sqrt{(\tilde{v}^{(1)}, \tilde{u}^{(1)})}$
\State for $j=1,2,\ldots,$ until convergence do
\State $\tilde{v}^{(j)}=\frac{\tilde{v}^{(j)}}{\tilde{\gamma}_j}$, $\tilde{u}^{(j)}=\frac{\tilde{u}^{(j)}}{\tilde{\gamma}_j}$
\State $\tilde{\delta}_j=(\tilde{u}^{(j)}, {\cal P}_i^{-1}{\cal A}_i\tilde{v}^{(j)})$
\State $\tilde{v}^{(j+1)}={\cal P}_i^{-1}{\cal A}_i\tilde{v}^{(j)}-\tilde{\delta}_j\tilde{v}^{(j)}-\tilde{\gamma}_j\tilde{v}^{(j-1)}$
\State $\tilde{u}^{(j+1)}={\cal P}_i\tilde{v}^{(j+1)}$
\State $\tilde{\gamma}_{j+1}=\sqrt{(\tilde{v}^{(j+1)}, \tilde{u}^{(j+1)})}$. If $\tilde{\gamma}_{j+1}=0$, go to line 13
\State $\tilde{T}_j=\mathrm{tridiag}(\tilde{\gamma}_{j}, \tilde{\delta}_j, \tilde{\gamma}_{j+1})\in R^{(j+1)\times j}$ and $e_1=(1,0,\ldots,0)^T\in R^{j+1}$
\State Calculate $x_j=x_0+V_jy^{(j)}$, where $V_j=[\tilde{v}^{(1)}, \ldots, \tilde{v}^{(j)}]\in R^{N\times j}$ and $y^{(j)}=\arg\min_{y\in R^{j}}\|\tilde{\gamma}_1e_1-\tilde{T}_jy\|_2$.
\State end do
\State $l:=j$
\State Calculate $x_l=x_0+V_ly^{(l)}$, where $V_l=[\tilde{v}^{(1)}, \ldots, \tilde{v}^{(l)}]\in R^{N\times l}$ and $y^{(l)}=\arg\min_{y\in R^{l}}\|\tilde{\gamma}_1e_1-\tilde{T}_ly\|_2$.
\end{algorithmic}
\end{algorithm}

Next, we are ready to validate the effectiveness of Algorithm 1 and Algorithm 2 by the sharp decaying numbers of iterations.

\subsection{Case 1: state constraint and lower control bound}\label{sec:case1}

We select $f, y_d, \tilde{g}_d$ as
 \begin{eqnarray*}
\left\{
\begin{aligned}
&f=(8\pi^2+\frac{1}{\alpha})\sin(2\pi x_1)\sin(2\pi x_2)-1,\\
&y_d = (1-\frac{8\pi^2}{\beta})\sin(2\pi x_1)\sin(2\pi x_2)-\frac{0.5}{\beta}, \\
&\tilde{g}_d = (2\pi\cos(2\pi x_1)\sin(2\pi x_2), 2\pi\sin(2\pi x_1)\cos(2\pi x_2))^T ,
\end{aligned}
\right.
\end{eqnarray*}
and set $\kappa=1$, $\gamma = 1$, $\beta = 1$.
Then the corresponding solution reads
\begin{eqnarray*}
\left\{
\begin{aligned}
&u=-\frac{1}{\alpha}\sin(2\pi x_1)\sin(2\pi x_2)+1,\\
&y = \sin(2\pi x_1)\sin(2\pi x_2), \\
&p = \sin(2\pi x_1)\sin(2\pi x_2), \\
&\lambda = -0.5,\\
&z = -\alpha\le0.
\end{aligned}
\right.
\end{eqnarray*}
Here $\int_\Omega u\,dx=1=\hat a$ and $z=-\alpha<0$ (the latter because $\alpha u+p=\alpha$), so that
the state constraint and the \emph{lower} control bound are both active and the complementarity
conditions are satisfied with the multipliers $\lambda=-1/2$ and $z=-\alpha$; this example therefore
exercises case $S5$ of Algorithm 1.

In this case, the elliptic optimal control problems are discretized by linear finite elements and
calculated with our proposed algorithms. The numerical results about relative error for the control,
state and co-state are depicted in Table~\ref{tab:re1} with $\alpha=1.0$e-6, where
\begin{eqnarray}\label{def_re}
u_{rel}=\frac{\|u-u_h\|_{L^2(\Omega)}}{\|u\|_{L^2(\Omega)}},\;\;\;y_{rel}=\frac{\|y-y_h\|_{H^1(\Omega)}}{\|y\|_{H^1(\Omega)}},\;\;\;p_{rel}=\frac{\|p-p_h\|_{H^1(\Omega)}}{\|p\|_{H^1(\Omega)}}.
\end{eqnarray}

\begin{table}[H]
\centering
\caption{Relative errors of control, state and co-state (Case~1) with $\alpha=1.0$e-6 and different $h$}
\label{tab:re1}
\begin{tabular}{c|cccc}
	\toprule
$\alpha$  &$h$            &  $u_{rel}$      & $y_{rel}$    & $p_{rel}$  \\ \hline
\multirow{4}{*}{$1.0$e-6}
&$6.250$e-2   &  1.182e-2         &   2.907e-1      &  3.235e-1   \\ \cmidrule(l){2-5}
&$3.125$e-2   &   2.634e-3        &    7.412e-2     &  7.437e-2   \\  \cmidrule(l){2-5}
&$1.563$e-2   &  6.273e-4  &    3.653e-2     &  3.653e-2   \\  \cmidrule(l){2-5}
&$7.813$e-3   &  1.532e-4  &    1.822e-2     &   1.822e-2   \\
	\bottomrule
\end{tabular}
\end{table}

By the data in Table~\ref{tab:re1} we conclude that the relative errors of the control, state and
co-state decrease as $h$ decreases. The observed orders of convergence are consistent with the use of
$P_1$ finite elements: $u_{rel}=O(h^2)$ in the $L^2$-norm, while $y_{rel},p_{rel}=O(h)$ in the
$H^1$-norm. These results show that the proposed algorithms solve this distributed convex optimal
control problem with integral constraints on the control and state effectively.

Furthermore, to assess the performance of our preconditioning strategy, we solve the algebraic systems (\ref{MATRix-1})-(\ref{MATRix-6}) using MINRES without preconditioning (MINRES/PRE, for short), MINRES with our exact preconditioner ${\cal P}_i$ (MINRES with ${\cal P}$, for short), and MINRES with our inexact preconditioner $\tilde{\cal P}_i$ (MINRES with $\tilde{\cal P}$, for short). We further include the block preconditioner from \cite{RDW} as a demanding and informative benchmark; this preconditioner is state-of-the-art for distributed PDE-constrained optimization problems without control constraints and forms a natural reference for our doubly constrained setting. For the distributed control system  
\[ \begin{pmatrix} 2\beta M & 0 & K^T\\ 0 & M & -M\\ K & -M & 0\end{pmatrix}, \]
\cite[eq.~(3.1)]{RDW} proposes the block-diagonal preconditioner
\[{\cal P}_{D1}=\begin{pmatrix} 2\beta M & 0 & 0\\ 0 & M & 0\\ 0 & 0 & KM^{-1}K^T\end{pmatrix}.\]
We transplant this construction to our system ${\cal A}_i$ as
\[
\mathcal P_1^{\rm RDW}=
\begin{pmatrix}
\alpha B & 0 & 0\\[2pt]
0 & S & 0\\[2pt]
0 & 0 & \mathsf C
\end{pmatrix},\quad \mathcal P_2^{\rm RDW}=
\begin{pmatrix}
\alpha B & 0 & 0 & 0\\[2pt]
0 & S & 0 & 0\\[2pt]
0 & 0 & \mathsf C & -\alpha^{-1} BD^{-1}L\\[2pt]
0 & 0 & -\alpha^{-1}L^TD^{-1} B & \alpha^{-1}L^TD^{-1}L
\end{pmatrix},
\]
\[
\mathcal P_3^{\rm RDW}=
\begin{pmatrix}
\alpha B & 0 & 0 & 0\\[2pt]
0 & S & 0 & 0\\[2pt]
0 & 0 & \mathsf C & A S^{-1}L\\[2pt]
0 & 0 & L^TS^{-1}A & L^TS^{-1}L
\end{pmatrix}, \quad \mathcal P_4^{\rm RDW}=
\begin{pmatrix}
\alpha B & 0 & 0 & 0 & 0\\[2pt]
0 & S & 0 & 0 & 0\\[2pt]
0 & 0 & \mathsf C & A S^{-1}L & -\alpha^{-1}B D^{-1}L\\[2pt]
0 & 0 & L^TS^{-1}A & L^TS^{-1}L & 0\\[2pt]
0 & 0 & -\alpha^{-1}L^TD^{-1} B & 0 & \alpha^{-1}L^TD^{-1}L
\end{pmatrix},
\]
where $S=\beta B+\gamma A$, and $D=\mathrm{diag}(B)$, $\mathsf C=\alpha^{-1} B D^{-1} B+A S^{-1}A$.
And we solve the algebraic systems (\ref{MATRix-1})-(\ref{MATRix-6}) using MINRES with these (MINRES with $\mathcal P^{\rm RDW}$, for short). All the iteration processes are started at the zero
vector and terminated if the current iterations satisfy tolerance $<1.0$e-6 or the maximum iteration
$> 5000$. By Algorithm~1, we know that we only need to calculate the first five steps. The numerical
results about iteration numbers, slove time and final relative residual for the systems (\ref{MATRix-1})-(\ref{MATRix-5}) (denoted by $S1, S2, S3, S4, S5$) with
different $\alpha$ and $h$ are shown in Table~\ref{tab:it1} and Table~\ref{tab:it3}.

\begin{table}[H]
\centering
\caption{Numbers of iterations for systems $S1$--$S5$ with different meshes and $\alpha=1.0$e-6 and $1.0$e-8}\label{tab:it1}
\resizebox{\textwidth}{!}{%
\begin{tabular}{cccccccccccc}
	\toprule
	\multirow{3}{*}{\textbf{$h$}} &\multirow{3}{*}{\shortstack{\text{Scheme}} } &\multicolumn{5}{c}{\textbf{$\alpha = 1.0$e-6}} &\multicolumn{5}{c}{\textbf{$\alpha = 1.0$e-8}} \cr
	\cmidrule(l){3-7} \cmidrule(l){8-12}
	&& $S1$ &$S2$ &$S3$
	&$S4$ &$S5$ & $S1$ &$S2$ &$S3$
	&$S4$ &$S5$\cr
	\midrule
	\multirow{4}{*}{6.250e-2}
	&MINRES/PRE. & 1027  & 973 & 987 & 1121  & 1279  & 942 &983 & 992 & 1120 & 1200\cr
	& MINRES with $\mathcal P^{\rm RDW}$ & 20  & 20 & 20 & 20  & 20  & 20 & 20 & 20 & 20 & 20\cr
	&MINRES with ${\cal P}$ & 3  & 3 & 3 & 3  & 3  & 2 & 2  & 2 & 2& 2 \cr
	&MINRES with $\tilde{\cal P}$ & 17  & 17 & 17 & 17  & 17  & 10 & 10  & 10 & 12& 12 \cr
\midrule
\multirow{4}{*}{3.125e-2}
	&MINRES/PRE. & 2491  & 2840 & 2844 & 2998  & 3327  & 2493 &2810 & 2849 & 3027 & 3351\cr
	&MINRES with $\mathcal P^{\rm RDW}$ & 18  & 18 & 18 & 18  & 18  & 18 & 18  & 18 & 18 & 18 \cr
	&MINRES with ${\cal P}$ & 3  & 3 & 3 & 3  & 3  & 2 & 2  & 2 & 2& 2  \cr
	&MINRES with $\tilde{\cal P}$ & 15  & 15 & 15 & 15  & 15  & 10 & 10  & 10 & 12& 12 \cr
\midrule
\multirow{4}{*}{1.563e-2}
	&MINRES/PRE. & 5000  & 5000 & 5000 & 5000  & 5000  & 5000  & 5000 & 5000 & 5000  & 5000\cr
	&MINRES with $\mathcal P^{\rm RDW}$ & 18  & 18 & 18 & 18  & 18  & 16 & 16  & 16 & 16& 16 \cr
	&MINRES with ${\cal P}$ & 3  & 3 & 3 & 3  & 3  & 2 & 2  & 2 & 2& 2 \cr
	&MINRES with $\tilde{\cal P}$ & 15  & 15 & 15 & 15  & 15  & 2 & 2  & 2 & 2& 2 \cr
	\bottomrule
\end{tabular}
}
\end{table}

\begin{table}[H]
\centering
\caption{Solve time (s) for systems $S1$--$S5$ with different meshes and $\alpha=1.0$e-6 and $1.0$e-8}\label{tab:it2}
\resizebox{\textwidth}{!}{%
\begin{tabular}{cccccccccccc}
	\toprule
	\multirow{3}{*}{\textbf{$h$}} &\multirow{3}{*}{\shortstack{\text{Scheme}} } &\multicolumn{5}{c}{\textbf{$\alpha = 1.0$e-6}} &\multicolumn{5}{c}{\textbf{$\alpha = 1.0$e-8}} \cr
	\cmidrule(l){3-7} \cmidrule(l){8-12}
	&& $S1$ &$S2$ &$S3$
	&$S4$ &$S5$ & $S1$ &$S2$ &$S3$
	&$S4$ &$S5$\cr
	\midrule
	\multirow{4}{*}{6.250e-2}
	&MINRES/PRE. & 0.0989  &  0.0115 & 0.00966  & 0.0111  & 0.0122 &0.0116 &   0.00891 &   0.00902 &    0.0115 &  0.0166 \cr
	& MINRES with $\mathcal P^{\rm RDW}$ & 0.0166  & 0.0023 & 0.00188 & 0.00172  & 0.00202  & 0.00207 &0.00197 &   0.00201 & 0.00199 & 0.00275\cr
	&MINRES with ${\cal P}$ & 0.00928  & 0.000608 & 0.00102 & 0.000689  & 0.000701  & 0.000648 & 0.000593 &  0.000715 & 0.00124 & 0.000702 \cr
	&MINRES with $\tilde{\cal P}$ & 0.00656 & 0.00149 & 0.00119 & 0.00132  & 0.00174  & 0.00115 & 0.00108 & 0.000907 & 0.00177 & 0.00162 \cr
\midrule
\multirow{4}{*}{3.125e-2}
	&MINRES/PRE. & 0.0853  & 0.0924 & 0.0881 & 0.101  & 0.106  & 0.098 &0.133 &0.134 & 0.119 & 0.132\cr
	&MINRES with $\mathcal P^{\rm RDW}$ & 0.0245 &  0.026 &0.0253 &0.0246 &0.023  & 0.0274 &0.0213 &0.03 &0.0257 &0.0206 \cr
	&MINRES with ${\cal P}$ & 0.00576 &0.00649 &0.00562 &0.00603 &0.00545& 0.00375 &0.00587 &0.00643 &0.00577 &0.00404  \cr
	&MINRES with $\tilde{\cal P}$ & 0.023 &0.0237 &0.0234 &0.0232 &0.0211  & 0.0111 &0.0185 &0.0186 &0.0135 &0.0133 \cr
\midrule
\multirow{4}{*}{1.563e-2}
	&MINRES/PRE. & 0.855 &      1.21 &      1.19 &       1.1 &      1.08  & 1.12 &1.06 &1.09 &0.993 &1.04\cr
	&MINRES with $\mathcal P^{\rm RDW}$ & 0.392 &     0.378 &     0.331 &    0.315  &  0.43  & 0.313 &0.281 &0.315 &0.277 &0.291 \cr
	&MINRES with ${\cal P}$ & 0.0893 &    0.0753 &     0.094 &    0.0738 &      0.08  & 0.0536 &0.0568 &0.051 &0.0645 &0.0534 \cr
	&MINRES with $\tilde{\cal P}$ & 0.308 & 0.329 & 0.312 &  0.287 &     0.306
  & 0.0483 &0.0605 &0.0519 &0.0473 &0.0584 \cr
	\bottomrule
\end{tabular}
}
\end{table}

\begin{table}[H]
\centering
\caption{Final relative residual for systems $S1$--$S5$ with different meshes and $\alpha=1.0$e-6 and $1.0$e-8}\label{tab:it3}
\resizebox{\textwidth}{!}{%
\begin{tabular}{cccccccccccc}
	\toprule
	\multirow{3}{*}{\textbf{$h$}} &\multirow{3}{*}{\shortstack{\text{Scheme}} } &\multicolumn{5}{c}{\textbf{$\alpha = 1.0$e-6}} &\multicolumn{5}{c}{\textbf{$\alpha = 1.0$e-8}} \cr
	\cmidrule(l){3-7} \cmidrule(l){8-12}
	&& $S1$ &$S2$ &$S3$
	&$S4$ &$S5$ & $S1$ &$S2$ &$S3$
	&$S4$ &$S5$\cr
	\midrule
	\multirow{4}{*}{6.250e-2}
	&MINRES/PRE. & 9.97e-07 &9.97e-07 &9.96e-07 &9.96e-07 &9.98e-07  & 9.88e-07 &9.91e-07 &9.82e-07 &9.79e-07 &1e-06\cr
	& MINRES with $\mathcal P^{\rm RDW}$ & 7.65e-07 &7.47e-07 &7.47e-07 &7.65e-07 &7.47e-07  & 7.94e-07 &7.76e-07 &7.76e-07 &7.94e-07 &7.76e-07\cr
	&MINRES with ${\cal P}$ & 8.36e-16 &7.3e-16 &1.58e-15 &8.21e-16 &1.15e-15  & 4.02e-08 &4.02e-08 &4.02e-08 &3.89e-08 &3.89e-08 \cr
	&MINRES with $\tilde{\cal P}$ & 5.93e-07 &5.69e-07 &5.69e-07 &3.28e-07 &3.36e-07  & 7.75e-07 &7.43e-07 &7.43e-07 &6.7e-07 &7.05e-07 \cr
\midrule
\multirow{4}{*}{3.125e-2}
	&MINRES/PRE. & 9.97e-07 &9.88e-07 &9.93e-07 &9.93e-07 &9.95e-07  & 9.97e-07 &9.97e-07 &9.88e-07 &9.79e-07 &9.93e-07\cr
	&MINRES with $\mathcal P^{\rm RDW}$ & 9.24e-07  &9.41e-07  & 9.39e-07 &9.85e-07 &9.91e-07  & 7.68e-07 &7.74e-07 &7.74e-07 &7.68e-07 &7.74e-07 \cr
	&MINRES with ${\cal P}$ & 1.5e-15 &2.76e-15 &5.22e-15 &1.62e-15 &1.55e-15  & 1.97e-08 &1.97e-08 &1.97e-08 &1.71e-08 &1.71e-08  \cr
	&MINRES with $\tilde{\cal P}$ & 6.85e-07 & 6.91e-07  &6.39e-07 &6.97e-07 &6.97e-07  & 7.49e-07 &7.47e-07 &7.47e-07 &4.95e-07 &5.03e-07 \cr
\midrule
\multirow{4}{*}{1.563e-2}
	&MINRES/PRE. & 3.35e-06 &1.51e-06 &1.61e-06 &3.9e-06 &2.29e-06  & 3.35e-06 &1.36e-06 &1.35e-06 &3.9e-06 &2.29e-06\cr
	&MINRES with $\mathcal P^{\rm RDW}$ & 6.95e-07 &6.96e-07 &6.91e-07 &6.82e-07 &6.84e-07  & 8.58e-07 &8.58e-07 &8.58e-07 &8.58e-07 &8.58e-07 \cr
	&MINRES with ${\cal P}$ & 2.73e-15 &2.79e-15 &6.87e-15 &2.75e-15 &5.31e-15  & 1.53e-08 &1.53e-08 &1.53e-08 &1.17e-08 &1.17e-08 \cr
	&MINRES with $\tilde{\cal P}$ & 4.9e-07 &4.92e-07 &4.9e-07 &4.91e-07 &4.92e-07  & 7.29e-07 &7.29e-07 &7.29e-07 &7.29e-07 &7.29e-07 \cr
	\bottomrule
\end{tabular}
}
\end{table}
Tables \ref{tab:it1}--\ref{tab:it3} present the iteration numbers , solve time and final relative residuals of four MINRES implementations for five linear systems $S1$--$S5$ with different mesh sizes $h$ and parameter $\alpha$. The unpreconditioned MINRES performs worst: its iteration count rises rapidly upon mesh refinement, hitting the maximum limit of 5000 for $h=1.563\times10^{-2}$ with considerably increased computational cost. The reference preconditioner $\mathcal P^{\rm RDW}$ from the literature effectively reduces iterations and exhibits good mesh robustness, yet it still requires more iterations than our exact preconditioner $\mathcal P$. The proposed exact preconditioner $\mathcal P$ achieves optimal performance across all test cases, with a stable iteration number of 2--3 that is nearly insensitive to mesh resolution and $\alpha$, the shortest solve time, and residuals close to machine precision. Our inexact preconditioner $\tilde{\mathcal P}$, constructed by incomplete Cholesky factorization, also outperforms unpreconditioned MINRES and $\mathcal P^{\rm RDW}$. For fine meshes with $\alpha=10^{-8}$, its iteration count matches that of the exact preconditioner, while enjoying lower assembly and factorization overhead to attain a favorable balance between efficiency and construction cost. Residual comparisons show that all preconditioned methods keep relative residuals within a reasonable range; $\tilde{\mathcal P}$ maintains strong acceleration with only mild residual degradation. Numerical results verify that both the exact and inexact variants of our preconditioning framework outperform the competitors and are suitable for solving this family of parameter-dependent linear systems.

\subsection{Case 2: state constraint and upper control bound}\label{sec:case2}

The given functions $f, y_d, \tilde{g}_d$ are selected as
 \begin{eqnarray*}
\left\{
\begin{aligned}
&f=8\pi^2\sin(2\pi x_1)\sin(2\pi x_2)-\frac{10}{\alpha}\sin(\pi x_1)\sin(\pi x_2)+\frac{40}{\alpha \pi^2}-2,\\
&y_d = \sin(2\pi x_1)\sin(2\pi x_2)+\frac{20\pi^2}{\beta}\sin(\pi x_1)\sin(\pi x_2)-\frac{1}{\beta}, \\
&\tilde{g}_d = (2\pi\cos(2\pi x_1)\sin(2\pi x_2), 2\pi\sin(2\pi x_1)\cos(2\pi x_2))^T ,
\end{aligned}
\right.
\end{eqnarray*}
and $\kappa=1$, $\gamma = 1$, $\beta = \frac{1}{2}$. Then the solution for this case reads
\begin{eqnarray*}
\left\{
\begin{aligned}
&u=\frac{10}{\alpha}\sin(\pi x_1)\sin(\pi x_2)-\frac{40}{\alpha \pi^2}+2,\\
&y = \sin(2\pi x_1)\sin(2\pi x_2), \\
&p = -10\sin(\pi x_1)\sin(\pi x_2), \\
&\lambda = -1,\\
&z = \frac{40}{\pi^2}-2\alpha>0.
\end{aligned}
\right.
\end{eqnarray*}
Here $\int_\Omega u\,dx=2=\hat b$ and $z=40/\pi^2-2\alpha>0$; consequently the state constraint and the
\emph{upper} control bound are both active, the multipliers are $\lambda=-1$ and
$z=40/\pi^2-2\alpha$, and this example exercises case $S6$ of Algorithm~1.

In this case, the constrained optimal control problems are calculated with our proposed Algorithm 1 and
Algorithm 2. The relative errors of the control, state and co-state variables are given in
Table~\ref{tab:re2} with $\alpha=1.0$e-6.

\begin{table}[H]
\centering
\caption{Relative errors of control, state and co-state (Case~2) with $\alpha=1.0$e-6 and different $h$}
\label{tab:re2}
\begin{tabular}{c|cccc}
	\toprule
$\alpha$  &$h$            &  $u_{rel}$      & $y_{rel}$    & $p_{rel}$  \\ \hline
\multirow{4}{*}{$1.0$e-6}
&$6.250$e-2   &  1.288e-1  &   1.451e-1      &  9.672e-2   \\ \cmidrule(l){2-5}
&$3.125$e-2   &  3.224e-2  & 7.274e-2        & 5.451e-2    \\  \cmidrule(l){2-5}
&$1.563$e-2   & 8.323e-3   &     3.643e-2    &  2.852e-2   \\  \cmidrule(l){2-5}
&$7.813$e-3   & 2.078e-3   &     1.824e-2    & 1.423e-2     \\
	\bottomrule
\end{tabular}
\end{table}

In view of the data in Table~\ref{tab:re2}, we find that the errors of the control, state and co-state
decrease with the decreased $h$, which illustrate the effectiveness of Algorithm 1 and Algorithm 2 for
solving the control-state constrained optimal control problem. We also investigate the efficiency of our
proposed preconditioners ${\cal P}_i$. Similarly, we calculate the algebraic systems (2.18a)--(2.18f) by
MINRES with preconditioners ${\cal P}_i$, $\tilde{\cal P}_i$  $(i=1,2,3,4)$ comparing with the MINRES method without
preconditioning and MINRES method with preconditioners $\mathcal P_i^{\rm RDW}$ $(i=1,2,3,4)$. The initial iteration value is zero vector and the iteration termination condition is
that either the tolerance $<1.0$e-6 or the maximum iteration $> 5000$. By Algorithm~1, we need to
iterate all 6 steps. The numerical
results about iteration numbers, slove time and final relative residual for the systems (\ref{MATRix-1})-(\ref{MATRix-6}) (denoted by $S1, S2, S3, S4, S5, S6$) with
different $\alpha$ and $h$ are shown in Table~\ref{tab:it4} and Table~\ref{tab:it6}.

\begin{table}[H]
\centering
\setlength{\tabcolsep}{3pt}
\caption{Numbers of iterations for systems $S1$--$S6$ with different meshes and $\alpha=1.0$e-6 and $1.0$e-8}\label{tab:it4}
\resizebox{\textwidth}{!}{%
\begin{tabular}{cccccccccccccc}
	\toprule
	\multirow{3}{*}{\textbf{$h$}} &\multirow{3}{*}{\shortstack{\text{Scheme}} } &\multicolumn{6}{c}{\textbf{$\alpha = 1.0$e-6}} &\multicolumn{6}{c}{\textbf{$\alpha = 1.0$e-8}} \cr
	\cmidrule(l){3-8} \cmidrule(l){9-14}
	&& $S1$ &$S2$ &$S3$
	&$S4$ &$S5$ & $S6$ & $S1$ &$S2$ &$S3$
	&$S4$ &$S5$ & $S6$ \cr
	\midrule
	\multirow{4}{*}{6.250e-2}
	&MINRES/PRE. & 2170  & 2557 & 2556 & 2411  & 2919  & 2922 &2051 &2445 &2417 &2334 & 2687 &  2683\cr
	&MINRES with $\mathcal P^{\rm RDW}$ & 22 &20 &20 &22 &20 &20  & 22 &20 &20 &22 &20 &20  \cr
	&MINRES with ${\cal P}$ & 3  & 3 & 3 & 3  & 3  & 3 & 2  & 2 & 2 & 2  & 2  & 2  \cr
	&MINRES with $\tilde{\cal P}$ & 18 &15 &15 &15 &13 &13 & 17 &13 &13 &14 &11 &11  \cr
\midrule
\multirow{4}{*}{3.125e-2}
	&MINRES/PRE. & 5000  & 5000 & 5000 & 5000  & 5000  & 5000  & 5000 & 5000 & 5000  & 5000 & 5000 & 5000 \cr
	&MINRES with $\mathcal P^{\rm RDW}$ & 29 &23 &23 &23 &18 &18  & 24 &18 &18 &22 &18 &18  \cr
	&MINRES with ${\cal P}$ & 3  & 3 & 3 & 3  & 3  & 3 & 2  & 2 & 2 & 2  & 2  & 2  \cr
	&MINRES with $\tilde{\cal P}$ & 22 &18 &18 &18 &16 &16 & 22 &10 &10 &16 &12 &12  \cr
\midrule
\multirow{4}{*}{1.563e-2}
	&MINRES/PRE. & 5000  & 5000 & 5000 & 5000  & 5000  & 5000  & 5000 & 5000 & 5000  & 5000 & 5000 & 5000 \cr
	&MINRES with $\mathcal P^{\rm RDW}$ & 50 &24 &24 &37 &22 &22  & 34 &16 &16 &27 &18 &18   \cr
	&MINRES with ${\cal P}$ & 3  & 3 & 3 & 3  & 3  & 3 & 2  & 2 & 2 & 2  & 2  & 2  \cr
	&MINRES with $\tilde{\cal P}$ & 20 &16 &16 &20 &16 &16 & 22 &15 &15 &22 &16 &16  \cr
	\bottomrule
\end{tabular}
}
\end{table}

\begin{table}[H]
\centering
\setlength{\tabcolsep}{3pt}
\caption{Solve time (s) for systems $S1$--$S6$ with different meshes and $\alpha=1.0$e-6 and $1.0$e-8}\label{tab:it5}
\resizebox{\textwidth}{!}{%
\begin{tabular}{cccccccccccccc}
	\toprule
	\multirow{3}{*}{\textbf{$h$}} &\multirow{3}{*}{\shortstack{\text{Scheme}} } &\multicolumn{6}{c}{\textbf{$\alpha = 1.0$e-6}} &\multicolumn{6}{c}{\textbf{$\alpha = 1.0$e-8}} \cr
	\cmidrule(l){3-8} \cmidrule(l){9-14}
	&& $S1$ &$S2$ &$S3$
	&$S4$ &$S5$ & $S6$ & $S1$ &$S2$ &$S3$
	&$S4$ &$S5$ & $S6$ \cr
	\midrule
	\multirow{4}{*}{6.250e-2}
	&MINRES/PRE. &  0.0146 &0.0174 &0.0187 &0.0171 &0.0204 &0.0215 & 0.0152 &0.0217 &0.0182 &0.0192 &0.0213 &0.0205 \cr
	&MINRES with $\mathcal P^{\rm RDW}$ & 0.00172 &0.00154 &0.00172 &0.00173 &0.00203 &0.00151 & 0.00209 &0.0021 &0.00272 &0.00168 &0.00158 &0.00184 \cr
	&MINRES with ${\cal P}$ & 0.000611 &0.000868 &0.000407 &0.0005 &0.000429 &0.000554 & 0.000697 &0.000472 &0.000773 &0.000463 &0.000628 &0.00049  \cr
	&MINRES with $\tilde{\cal P}$ & 0.0014 &0.00102 &0.000926 &0.00121 &0.000875 &0.000875 & 0.00156 &0.000975 &0.0013 &0.00113 &0.00101 &0.000864   \cr
\midrule
\multirow{4}{*}{3.125e-2}
	&MINRES/PRE. & 0.151 &0.149 &0.147 &0.148 &0.156 &0.153 & 0.152 &0.155 &0.149 &0.148 &0.157 &0.156  \cr
	&MINRES with $\mathcal P^{\rm RDW}$ & 0.0275 &0.0256 &0.0275 &0.025 &0.0198 &0.022 & 0.0385 &0.0215 &0.0203 &0.0237 &0.0222 &0.0264   \cr
	&MINRES with ${\cal P}$ & 0.00511 &0.00475 &0.0065 &0.00707 &0.00481 &0.00558 & 0.00497 &0.00395 &0.00387 &0.00454 &0.0051 &0.00556  \cr
	&MINRES with $\tilde{\cal P}$ & 0.0243 &0.0205 &0.0232 &0.022 &0.0162 &0.0151  & 0.0262 & 0.0121 &0.0112 &0.0158 &0.0113 &0.0139    \cr
\midrule
\multirow{4}{*}{1.563e-2}
	&MINRES/PRE. & 1.02 &1.15 &1.01 &0.971 &1.03 &1.04 & 1.03 &0.967 &0.987 &1.25 &1.31 &1.51  \cr
	&MINRES with $\mathcal P^{\rm RDW}$ & 0.818 &0.377 &0.367 &0.592 &0.394 &0.341 & 0.514 &0.286 &0.562 &0.839 &0.548 &0.484    \cr
	&MINRES with ${\cal P}$ & 0.0663 &0.0678 &0.0744 &0.0654 &0.0707 &0.074 & 0.0545 &0.0516 &0.114 &0.0772 &0.0939 &0.0885  \cr
	&MINRES with $\tilde{\cal P}$ & 0.341 &0.281 &0.276 &0.342 &0.291 &0.278 & 0.294 & 0.238 &0.463 &0.527 &0.489 &0.442    \cr
	\bottomrule
\end{tabular}
}
\end{table}

\begin{table}[H]
\centering
\setlength{\tabcolsep}{3pt}
\caption{Final relative residual for systems $S1$--$S6$ with different meshes and $\alpha=1.0$e-6 and $1.0$e-8}\label{tab:it6}
\resizebox{\textwidth}{!}{%
\begin{tabular}{cccccccccccccc}
	\toprule
	\multirow{3}{*}{\textbf{$h$}} &\multirow{3}{*}{\shortstack{\text{Scheme}} } &\multicolumn{6}{c}{\textbf{$\alpha = 1.0$e-6}} &\multicolumn{6}{c}{\textbf{$\alpha = 1.0$e-8}} \cr
	\cmidrule(l){3-8} \cmidrule(l){9-14}
	&& $S1$ &$S2$ &$S3$
	&$S4$ &$S5$ & $S6$ & $S1$ &$S2$ &$S3$
	&$S4$ &$S5$ & $S6$ \cr
	\midrule
	\multirow{4}{*}{6.250e-2}
	&MINRES/PRE. & 9.57e-07 &9.99e-07 &9.98e-07 &9.94e-07 &9.96e-07 &9.55e-07 & 9.99e-07 &  9.75e-07 &  9.84e-07 &  9.87e-07 &  9.93e-07&  9.98e-07  \cr
	&MINRES with $\mathcal P^{\rm RDW}$ & 9.24e-07 & 6.98e-07 &6.98e-07 &9.24e-07 &7.06e-07 &7.06e-07 & 9.68e-07 &  7.22e-07 &  7.22e-07 & 9.68e-07 &  7.22e-07&  7.22e-07   \cr
	&MINRES with ${\cal P}$ & 6.82e-16 &1.19e-15 &9.15e-16 &6.76e-16 &1.2e-15 &4.11e-15 & 4.07e-07 &  4.07e-07 &  4.07e-07 &  2.47e-07 &  2.47e-07 &  2.47e-07  \cr
	&MINRES with $\tilde{\cal P}$ & 4.04e-07 &4.46e-07 &4.49e-07 &7.91e-07 &  6.37e-07 &  6.24e-07 & 3.06e-07&  6.11e-07&  6.11e-07& 6.03e-07 &  6.21e-07& 6.21e-07   \cr
\midrule
\multirow{4}{*}{3.125e-2}
	&MINRES/PRE. & 1.08e-05 &  7.77e-06 &  8.49e-06 &   2.1e-05 &  7.39e-06 &  7.46e-06  & 1.08e-05& 8.52e-06 &  8.72e-06 &  1.74e-05 &  7.37e-06 & 7.33e-06  \cr
	&MINRES with $\mathcal P^{\rm RDW}$ &  8.95e-07&  9.74e-07 & 9.74e-07 &7.61e-07   & 8.38e-07& 8.49e-07 & 7.74e-07&  8.08e-07 &  8.08e-07 &  7.74e-07 &  8.13e-07 &  8.13e-07   \cr
	&MINRES with ${\cal P}$ &1.23e-15 & 9.09e-15 &  2.38e-15 &  1.53e-15 &  6.54e-15 & 1.29e-14 & 3.83e-07 &  3.83e-07&  3.83e-07&  2.56e-07 &  2.56e-07&  2.56e-07  \cr
	&MINRES with $\tilde{\cal P}$ & 6.64e-07 &  6.93e-07 & 6.93e-07 & 6.45e-07  &  7.18e-07&  7.09e-07 & 5.93e-07&  7.04e-07& 7.04e-07 & 7.58e-07 &   6.52e-07 &   6.53e-07   \cr
\midrule
\multirow{4}{*}{1.563e-2}
	&MINRES/PRE. &  0.000303&  0.000115& 9.75e-05&  0.000325&  3.25e-05& 3.66e-05 & 0.000303&  9.94e-05 &  0.000109 &  0.000325 &  3.62e-05 &  3.55e-05  \cr
	&MINRES with $\mathcal P^{\rm RDW}$ & 9.49e-07 & 8.05e-07 &  8.06e-07 &   9.5e-07&  8.35e-07 &  8.34e-07  & 4.89e-07 &  8.65e-07 &  8.65e-07 &  4.89e-07 &  8.85e-07& 8.85e-07  \cr
	&MINRES with ${\cal P}$ &  2.71e-15 &    4e-15 &  4.26e-14 &  2.68e-15 &   5.5e-14 & 3.81e-14 & 3.72e-07 &3.72e-07 &3.72e-07 &2.62e-07 &2.62e-07 &2.62e-07  \cr
	&MINRES with $\tilde{\cal P}$ & 7.66e-07 & 5.56e-07 &  5.55e-07&  4.34e-07&  7.1e-07 & 7.11e-07 & 3.31e-07 & 2.89e-07&  2.89e-07&  3.92e-07 &  4.59e-07&  4.59e-07   \cr
	\bottomrule
\end{tabular}
}
\end{table}
Tables \ref{tab:it4}-\ref{tab:it6} report the iteration number, solve time and final relative residuals for the systems $S1$--$S6$ with various mesh sizes $h$ and parameter $\alpha$. The numerical behavior is broadly consistent with case 1. The unpreconditioned MINRES quickly hits the maximum iteration limit under mesh refinement, while $\mathcal P^{\rm RDW}$ accelerates convergence but is inferior to our proposed preconditioners. The exact preconditioner $\mathcal P$ retains robust convergence with only 2--3 iterations, the lowest computing time and the best residual accuracy. The inexact preconditioner $\tilde{\mathcal P}$ also outperforms the competing schemes and achieves a good trade-off between convergence speed and construction cost. These results further confirm the effectiveness and robustness of our preconditioning framework for this class of linear systems.

\section{Concluding remarks}

In this work, we study eight valid preconditioners for the discretized algebraic systems generated from linear finite element methods of control-state constrained optimal control problems. Both rigorous theoretical analyses and numerical verifications are given to state our preconditioners are robust with respect to given parameters and mesh sizes.
Furthermore, the cost of iterations within the numerical tests is optimal, which affirms our proposed preconditioners significantly reduce the computational cost during iterations.
Meanwhile, our techniques can be extended to solve different constrained optimal control problems by other numerical methods, which are equipped  with the same saddle point systems.


\section*{Conflicts of Interest}
The authors declare no conflicts of interest.

\section*{Data Availability Statement}
All the data are available  within the article.


\end{document}